\documentclass{amsart}
\usepackage{amssymb,latexsym,hyperref}

\numberwithin{equation}{section}
\newtheorem{theorem}{Theorem}[section]

\newtheorem{remark}[theorem]{Remark}
\newtheorem{lemma}[theorem]{Lemma}

\begin{document}

\title
{Generalization of Schensted insertion algorithm to the cases of
hooks and semi-shuffles.}

\author {Mikhail Kogan}

\thanks{The author is supported by NSF Postdoctoral Fellowship.}

\address{\noindent Department of Mathematics, Northeastern University,
  Boston, MA 02115}
\email{misha@neu.edu}

\begin{abstract}
Given an rc-graph $R$ of permutation $w$ and an rc-graph $Y$ of
permutation $v$, we provide an insertion algorithm, which defines
an rc-graph $R\leftarrow Y$ in the case when $v$ is a shuffle with
the descent at $r$ and $w$ has no descents greater than $r$ or in
the case when $v$ is a shuffle, whose shape is a hook. This
algorithm gives a combinatorial rule for computing the generalized
Littlewood-Richardson coefficients~$c^{u}_{wv}$ in the two cases
mentioned above.
\end{abstract}

\maketitle

\pagestyle{myheadings} \markboth{MIKHAIL KOGAN}{GENERALIZATION OF
SCHENSTED INSERTION ALGORITHM.}


\section{Introduction.}
Rc-graphs were originally introduced by Fomin and Kirillov in
\cite{f-k}. They are explicit combinatorial objects, which encode
monomials in Schubert polynomials. Rc-graphs proved to be very
useful for providing combinatorial rules of computing certain
generalized Littlewood-Richardson (or just LR) coefficients (see
\cite{b-b}, \cite{k}, \cite{k-k}). In this paper we extend these
results to more general cases.

Denote by $\mathfrak S_w$ the Schubert polynomial of permutation
$w\in S_n$. Then the generalized LR coefficients $c^u_{wv}$ for
$u,v,w\in S_n$ are defined by
$$
\mathfrak S_w\cdot\mathfrak S_v=\sum_{u}c^u_{wv} \mathfrak S_u.
$$
(If $u,v,w$ are shuffles (also called grassmanian permutations)
with descents at $r$, the coefficients $c^{u}_{wv}$ are just the
LR coefficients.) It can be shown that all $c^u_{wv}$ are
nonnegative integers. (Consider the Schubert basis of the
cohomology ring of the flag variety, then $c^u_{wv}$ are the
structure constants, and they count the number of points in
certain intersections of algebraic varieties.) There are many
totally positive rules for computing LR coefficients (see \cite{f}
for further references), but there is no know totally positive
rule for generalized LR coefficients. (By a totally positive rule
we understand a construction of an explicit combinatorial set for
each triple $(u,v,w)$, such that $c^u_{wv}$ is equal to the number
of elements in this set.)

In certain cases (see \cite{b-s}, \cite{kn}) a totally positive
rule can by given by equating generalized LR coefficients to LR
coefficients. In other cases, such as Pieri formula (see
\cite{l-s}), a totally positive rule is given in terms of paths in
the Bruhat order.

Yet another approach to produce totally positive rule, adopted in
\cite{b-b}, \cite{k}, \cite{k-k} and in this paper, is to
generalize Schensted insertion algorithm to rc-graphs. The rule,
which we believe will be eventually generalized to the most
general case, is the following. An algorithm is constructed, which
inserts an rc-graph $Y$ of $v$ into rc-graph $R$ of $w$ to produce
an rc-graph $R\leftarrow Y$. Then for a fixed rc-graph $U$ of $u$,
$c^u_{wv}$ is the number of tuples $(R,Y)$ with $U=R\leftarrow Y$.
We present such an algorithm in the cases, when $v$ is an
$r$-shuffle and $w$ is an $r$-semi-shuffle, or when $v$ is a
shuffle, whose shape is a hook. (An $r$-shuffle is a shuffle with
the descent at $r$, an $r$-semi-shuffle is a permutation with no
descents greater than $r$.)

The first such algorithm was constructed by Bergeron and Billey
\cite{b-b} to prove Monk's formula (the case when $v$ is a simple
transposition). The author \cite{k} showed that their algorithm
works in the case when $w$ is an $r$-semi-shuffle and $v$ is an
$r$-shuffle. A modified algorithm was constructed by Kumar and the
author \cite{k-k} when $v$ is an $r$-shuffle whose shape is a row.
Kumar \cite{ku} constructed an analogous algorithm for the column
case. The analogy between two algorithms in \cite{k-k} and
\cite{ku} is similar to the analogy between row and column
insertion algorithms for Young tableaux.

This paper presents an algorithm, which works in all mentioned
above cases as well as in the case when $v$ is a shuffle, whose
shape is a hook. (A rule in this case written in terms of
$r$-Bruhat chains was originally constructed by Sottile \cite{s}.)
Our algorithm directly generalizes the algorithm of \cite{k-k}. In
the case when both $w$ and $v$ are $r$-shuffles, the algorithm
produces the same results as Schensted insertion algorithm. In all
the cases (see \cite{b-b}, \cite{k}, \cite{k-k}), except for the
case when $v$ is a shuffle whose shape is a~hook, the algorithm
can be simplified

Using the new insertion algorithm we also provide a rule for
computing generalized LR coefficients in the cases mentioned above
using $r$-Bruhat chains. This rule can be thought of as a
generalized RSK correspondence. In the case when the shape of $v$
is a hook, it is just a restatement of Pieri formula (see
\cite{s}). In the case when $v$ is an $r$-shuffle and $w$ is an
$r$-semi-shuffle it is a new result.

The paper is organized as follows. Section \ref{definitions}
introduces most of notations and definitions and contains the
statements of main results in
Theorems~\ref{thm:main},~\ref{thm:inverse
main}~and~\ref{thm:bruhat}. Theorem \ref{thm:main} states that the
algorithm defined in Section \ref{algorithm} works, it is proved
in Section \ref{proof}. Theorem \ref{thm:inverse main} states that
the inverse algorithm defined in Section~\ref{inverse algorithm}
works, it is given without a proof, since the proof is very
similar to the proof of Theorem~\ref{thm:main}. Theorem
\ref{thm:bruhat} gives a rule of computing certain generalized LR
coefficients in terms of $r$-Bruhat chains. Finally,
Section~\ref{examples} contains examples of the algorithm.
Section~\ref{examples} should be read together with Section
\ref{algorithm} to understand how the algorithm~works.

\smallskip

\noindent{\bf Acknowledgements.} The author thanks Frank Sotille,
Sara Billey and Abhinav Kumar for useful discussions as well as
Cristian Lenart for making valuable remarks about the paper.

\section{Notation, Definitions and Main Results.}
\label{definitions}
\subsection{Permutations}
\label{permutation} Let $S_n$ be the group of permutations
$w=(w(1),\dots,w(n))$ and let $S_\infty=\cup_nS_n$ be the group of
permutations on $\mathbb N$ which fix all but finitely many
integers. For $1\leq i<j$, denote by $t_{ij}$ the transposition,
which exchanges $i$ and $j$. The simple transpositions
$s_i=t_{i,i+1}$ for $1\leq i\leq n-1$ generate $S_n$.

A word $i_1\dots i_l$ in the alphabet $[1,2,\dots]$ is a reduced
word of $w\in S_\infty$, if $w=s_{i_1}\dots s_{i_l}$ and $l$ is
minimal. The length $l(w)$ of $w$ is set to be~$l$. The longest
permutation $w^n_0$ of the group $S_n$ is given by
$w^n_0(i)=n+1-i$ for $i\leq n$.

$w\in S_\infty$ is an \emph{$r$-shuffle} if $w(i)<w(i+1)$ for
$i\neq r$. It is an \emph{$r$-semi-shuffle} if $w(i)<w(i+1)$ for
$i>r$ . To each shuffle we associate a partition
$\lambda=(\lambda_1\geq \dots \geq \lambda_{r'}> 0)$ given by
$\lambda_j=w(r+1-j)-r-1+j$ for $j\leq r$, where
$\lambda_{r'+1}=0$. Then $w$ is uniquely determined by $\lambda$
and $r$ and we write $w=v(\lambda,r)$. A partition $\lambda$ is
\emph{a row} if $\lambda=(\lambda_1)$, it is \emph{a column} if
$\lambda=(1,\dots,1)$, and it is \emph{a hook} if
$\lambda=(\lambda_1,1,\dots,1)$.

\subsection{Rc-graphs.}
Let $W^{n}_0$ be the reduced word $(n-1\dots 1 \dots n-1\ n-2 \
n-1)$ of $w^{n}_0$. A subword $R$ of $W^n_0$ is called {\emph {a
graph}}. Each graph $R=i_1\dots i_m$ defines a permutation
$w(R)=s_{i_1}\dots s_{i_m}$. If $R$ is a reduced word of $w(R)$,
it is called {\it an rc-graph} of $w(R)$. Note that two different
subwords of $W^n_0$, which produce the same words are two
different graphs. For example, if $n=3$, $w=s_2$, then $W^3_0=212$
has two different subwords, whose permutation is $w$, namely the
subword $2$ placed at the first or third slot. The set of all
rc-graphs of $w$ is denoted by~$\mathcal R\mathcal C(w)$.

We think of graphs using the following pictorial presentation.
Think of $W^n_0$ as a triangular set of crossings shown in the
first picture of Figure \ref{definition of graphs} for $n=5$. Each
crossing is labelled by a letter from the alphabet
$[1,\dots,n-1]$. To get back $W^n_0$ we read those labels from top
to bottom row, from right to left in each row. Then each subword
$R$ of $W^n_0$ is presented as a subset of the crossings for
$W^n_0$. Two illustrations are provided in Figure \ref{definition
of graphs}, where the second picture corresponds to subword
$2323$, while the third picture corresponds to subword $4132$.

\begin{figure}[ht]

\begin{picture}(300,90)


\put(0,64){1} \put(0,48){2} \put(0,32){3} \put(0,16){4}
\put(0,0){5}

\put(18,83){1} \put(34,83){2} \put(50,83){3} \put(66,83){4}
\put(82,83){5}


\put (18,19){\line(1,0){6}} \put (21,16){\line(0,1){6}}
\put(22,14){\tiny{$4$}}

\put (18,35){\line(1,0){6}} \put (21,32){\line(0,1){6}}
\put(22,30){\tiny{$3$}}

\put (18,51){\line(1,0){6}} \put (21,48){\line(0,1){6}}
\put(22,46){\tiny{$2$}}

\put (18,67){\line(1,0){6}} \put (21,64){\line(0,1){6}}
\put(22,62){\tiny{$1$}}

\put (34,35){\line(1,0){6}} \put (37,32){\line(0,1){6}}
\put(38,30){\tiny{$4$}}

\put (34,51){\line(1,0){6}} \put (37,48){\line(0,1){6}}
\put(38,46){\tiny{$3$}}

\put (34,67){\line(1,0){6}} \put (37,64){\line(0,1){6}}
\put(38,62){\tiny{$2$}}

\put (50,51){\line(1,0){6}} \put (53,48){\line(0,1){6}}
\put(54,46){\tiny{$4$}}

\put (50,67){\line(1,0){6}} \put (53,64){\line(0,1){6}}
\put(54,62){\tiny{$3$}}

\put (66,67){\line(1,0){6}} \put (69,64){\line(0,1){6}}
\put(70,62){\tiny{$4$}}



\put(100,64){1} \put(100,48){2} \put(100,32){3} \put(100,16){4}
\put(100,0){5}

\put(118,83){1} \put(134,83){2} \put(150,83){3} \put(166,83){4}
\put(182,83){5}

\put (121,67){\circle*{3}} \put (121,19){\circle*{3}} \put
(137,35){\circle*{3}} \put (153,67){\circle*{3}} \put
(153,51){\circle*{3}}  \put (169,67){\circle*{3}}


\put (118,51){\line(1,0){6}} \put (121,48){\line(0,1){6}}
\put(122,46){\tiny{$2$}}

\put (118,35){\line(1,0){6}} \put (121,32){\line(0,1){6}}
\put(122,30){\tiny{$3$}}

\put (134,67){\line(1,0){6}} \put (137,64){\line(0,1){6}}
\put(138,62){\tiny{$2$}}

\put (134,51){\line(1,0){6}} \put (137,48){\line(0,1){6}}
\put(138,46){\tiny{$3$}}



\put(200,64){1} \put(200,48){2} \put(200,32){3} \put(200,16){4}
\put(200,0){5}

\put(218,83){1} \put(234,83){2} \put(250,83){3} \put(266,83){4}
\put(282,83){5}

\put (221,35){\circle*{3}} \put (221,19){\circle*{3}} \put
(237,35){\circle*{3}} \put (253,67){\circle*{3}} \put
(253,51){\circle*{3}}  \put (237,67){\circle*{3}}


\put (218,67){\line(1,0){6}} \put (221,64){\line(0,1){6}}
\put(222,62){\tiny{$1$}}

\put (218,51){\line(1,0){6}} \put (221,48){\line(0,1){6}}
\put(222,46){\tiny{$2$}}

\put (266,67){\line(1,0){6}} \put (269,64){\line(0,1){6}}
\put(270,62){\tiny{$4$}}

\put (234,51){\line(1,0){6}} \put (237,48){\line(0,1){6}}
\put(238,46){\tiny{$3$}}

\end{picture}
\caption{Examples of graphs.} \label{definition of graphs}
\end{figure}

Connect the crossings of $R$ by strands, which intersect at the
places, where there is a crossing and do not intersect otherwise.
For illustration see Figure \ref{more graphs}, where the graphs
correspond to the graphs from Figure \ref{definition of graphs}.
Notice that when we draw pictures of graphs we omit those parts of
graphs, which have no crossings. So graphs from Figures
\ref{definition of graphs} and \ref{more graphs} can be extended
to the right and down by nonintersecting strands.

\begin{figure}[ht]
\begin{picture}(300,90)

\put(0,64){1} \put(0,48){2} \put(0,32){3} \put(0,16){4}
\put(0,0){5}

\put(18,83){1} \put(34,83){2} \put(50,83){3} \put(66,83){4}
\put(82,83){5}

\put(81,71){\oval(8,8)[br]} \put(65,55){\oval(8,8)[br]}
\put(49,39){\oval(8,8)[br]} \put(33,23){\oval(8,8)[br]}
\put(17,7){\oval(8,8)[br]}


\put (18,19){\line(1,0){6}} \put (21,16){\line(0,1){6}}

\put (18,35){\line(1,0){6}} \put (21,32){\line(0,1){6}}

\put (18,51){\line(1,0){6}} \put (21,48){\line(0,1){6}}

\put (18,67){\line(1,0){6}} \put (21,64){\line(0,1){6}}

\put (34,35){\line(1,0){6}} \put (37,32){\line(0,1){6}}

\put (34,51){\line(1,0){6}} \put (37,48){\line(0,1){6}}

\put (34,67){\line(1,0){6}} \put (37,64){\line(0,1){6}}

\put (50,51){\line(1,0){6}} \put (53,48){\line(0,1){6}}

\put (50,67){\line(1,0){6}} \put (53,64){\line(0,1){6}}

\put (66,67){\line(1,0){6}} \put (69,64){\line(0,1){6}}


\put (08,67){\line(1,0){10}} \put (08,51){\line(1,0){10}} \put
(08,35){\line(1,0){10}} \put (08,19){\line(1,0){10}} \put
(08,3){\line(1,0){10}} \put (24,67){\line(1,0){10}} \put
(24,51){\line(1,0){10}} \put (24,35){\line(1,0){10}} \put
(24,19){\line(1,0){10}} \put (40,67){\line(1,0){10}} \put
(40,51){\line(1,0){10}} \put (40,35){\line(1,0){10}} \put
(56,67){\line(1,0){10}} \put (56,51){\line(1,0){10}} \put
(72,67){\line(1,0){10}}


\put (21,70){\line(0,1){10}} \put (21,54){\line(0,1){10}} \put
(21,38){\line(0,1){10}} \put (21,22){\line(0,1){10}} \put
(21,6){\line(0,1){10}} \put (37,70){\line(0,1){10}} \put
(37,54){\line(0,1){10}} \put (37,38){\line(0,1){10}} \put
(37,22){\line(0,1){10}} \put (53,70){\line(0,1){10}} \put
(53,54){\line(0,1){10}} \put (53,38){\line(0,1){10}} \put
(69,70){\line(0,1){10}} \put (69,54){\line(0,1){10}} \put
(85,70){\line(0,1){10}}



\put(100,64){1} \put(100,48){2} \put(100,32){3} \put(100,16){4}
\put(100,0){5}

\put(118,83){1} \put(134,83){2} \put(150,83){3} \put(166,83){4}
\put(182,83){5}

\put(181,71){\oval(8,8)[br]} \put(165,55){\oval(8,8)[br]}
\put(149,39){\oval(8,8)[br]} \put(133,23){\oval(8,8)[br]}
\put(117,7){\oval(8,8)[br]}

\put(165,71){\oval(8,8)[br]} \put(173,63){\oval(8,8)[tl]}
\put(149,71){\oval(8,8)[br]} \put(157,63){\oval(8,8)[tl]}
\put(149,55){\oval(8,8)[br]} \put(157,47){\oval(8,8)[tl]}
\put(133,39){\oval(8,8)[br]} \put(141,31){\oval(8,8)[tl]}
\put(117,23){\oval(8,8)[br]} \put(125,15){\oval(8,8)[tl]}
\put(117,71){\oval(8,8)[br]} \put(125,63){\oval(8,8)[tl]}


\put (118,51){\line(1,0){6}} \put (121,48){\line(0,1){6}}

\put (118,35){\line(1,0){6}} \put (121,32){\line(0,1){6}}

\put (134,67){\line(1,0){6}} \put (137,64){\line(0,1){6}}

\put (134,51){\line(1,0){6}} \put (137,48){\line(0,1){6}}


\put (108,67){\line(1,0){10}} \put (108,51){\line(1,0){10}} \put
(108,35){\line(1,0){10}} \put (108,19){\line(1,0){10}} \put
(108,3){\line(1,0){10}} \put (124,67){\line(1,0){10}} \put
(124,51){\line(1,0){10}} \put (124,35){\line(1,0){10}} \put
(124,19){\line(1,0){10}} \put (140,67){\line(1,0){10}} \put
(140,51){\line(1,0){10}} \put (140,35){\line(1,0){10}} \put
(156,67){\line(1,0){10}} \put (156,51){\line(1,0){10}} \put
(172,67){\line(1,0){10}}


\put (121,70){\line(0,1){10}} \put (121,54){\line(0,1){10}} \put
(121,38){\line(0,1){10}} \put (121,22){\line(0,1){10}} \put
(121,6){\line(0,1){10}} \put (137,70){\line(0,1){10}} \put
(137,54){\line(0,1){10}} \put (137,38){\line(0,1){10}} \put
(137,22){\line(0,1){10}} \put (153,70){\line(0,1){10}} \put
(153,54){\line(0,1){10}} \put (153,38){\line(0,1){10}} \put
(169,70){\line(0,1){10}} \put (169,54){\line(0,1){10}} \put
(185,70){\line(0,1){10}}



\put(200,64){1} \put(200,48){2} \put(200,32){3} \put(200,16){4}
\put(200,0){5}

\put(218,83){1} \put(234,83){2} \put(250,83){3} \put(266,83){4}
\put(282,83){5}

\put(281,71){\oval(8,8)[br]} \put(265,55){\oval(8,8)[br]}
\put(249,39){\oval(8,8)[br]} \put(233,23){\oval(8,8)[br]}
\put(217,7){\oval(8,8)[br]}

\put(233,71){\oval(8,8)[br]} \put(241,63){\oval(8,8)[tl]}
\put(249,71){\oval(8,8)[br]} \put(257,63){\oval(8,8)[tl]}
\put(249,55){\oval(8,8)[br]} \put(257,47){\oval(8,8)[tl]}
\put(233,39){\oval(8,8)[br]} \put(241,31){\oval(8,8)[tl]}
\put(217,23){\oval(8,8)[br]} \put(225,15){\oval(8,8)[tl]}
\put(217,39){\oval(8,8)[br]} \put(225,31){\oval(8,8)[tl]}


\put (218,67){\line(1,0){6}} \put (221,64){\line(0,1){6}}

\put (218,51){\line(1,0){6}} \put (221,48){\line(0,1){6}}

\put (266,67){\line(1,0){6}} \put (269,64){\line(0,1){6}}

\put (234,51){\line(1,0){6}} \put (237,48){\line(0,1){6}}


\put (208,67){\line(1,0){10}} \put (208,51){\line(1,0){10}} \put
(208,35){\line(1,0){10}} \put (208,19){\line(1,0){10}} \put
(208,3){\line(1,0){10}} \put (224,67){\line(1,0){10}} \put
(224,51){\line(1,0){10}} \put (224,35){\line(1,0){10}} \put
(224,19){\line(1,0){10}} \put (240,67){\line(1,0){10}} \put
(240,51){\line(1,0){10}} \put (240,35){\line(1,0){10}} \put
(256,67){\line(1,0){10}} \put (256,51){\line(1,0){10}} \put
(272,67){\line(1,0){10}}


\put (221,70){\line(0,1){10}} \put (221,54){\line(0,1){10}} \put
(221,38){\line(0,1){10}} \put (221,22){\line(0,1){10}} \put
(221,6){\line(0,1){10}} \put (237,70){\line(0,1){10}} \put
(237,54){\line(0,1){10}} \put (237,38){\line(0,1){10}} \put
(237,22){\line(0,1){10}} \put (253,70){\line(0,1){10}} \put
(253,54){\line(0,1){10}} \put (253,38){\line(0,1){10}} \put
(269,70){\line(0,1){10}} \put (269,54){\line(0,1){10}} \put
(285,70){\line(0,1){10}}
\end{picture}
\caption{Examples of graphs.} \label{more graphs}
\end{figure}

It is easy to see that for each graph $R$ and $i\in \mathbb N$,
$w(R)(i)$ is given by the column where the strand, which starts at
row $i$, ends. (Instead of referring to a strand as "a strand,
which starts at row $i$", we will say "strand $i$". So, the above
statement transforms to: strand $i$ ends in column $w(R)(i)$.)
This immediately leads to
\begin{align}
&\label{simple-add-crossing} \text{If graph $R$ is constructed
out of graph $R'$ by adding or removing \ \ \ \ \ \ \ \ \ \ \ \ } \\
&\nonumber\text{a crossing of strands $c$ and $d$ then
$w(R)=w(R')t_{cd}$.}
\end{align}

Denote by $|R|$ the length of the corresponding subword, or, in
other words, the number of crossings in $R$. Clearly, $|R|\geq
l(w(R))$. The following two statements are very easy to check and
are given without proofs.
\begin{align}
&\label{two crossings}\text{$R$ is an rc-graph if and only if
no two strands intersect twice in $R$.}\ \ \ \ \ \ \ \ \ \ \  \\
&\label{number-length} R\text{ is an rc-graph if and only if
}|R|=l(w(R)).
\end{align}
For example, the first and third graphs of Figure \ref{more
graphs} are rc-graphs, since both graphs have no double crossings
and for both $|R|=l(w(R))$. But the second graph is not an
rc-graph, since strands $3$ and $4$ intersect twice, or
$|R|=4>l(w(R))=2$.

Given a graph $R$, and a subset $I\subset \mathbb N$, define $R_I$
to be the graph, which coincides with $R$ at the rows labelled by
elements of $I$ and has no crossings outside these rows. For
example, if $I=\{\ell,\ell+1,\dots\}$, then $R_I= R_{\geq \ell}$
is the graph, which coincides with $R$ below or at row $\ell$ and
has no crossings above row $\ell$. Or, if $I=\{\ell\}$, then
$R_I=R_{\ell}$ coincides with $R$ at row $\ell$ and has no
crossings outside row $\ell$.

Given two graphs $R,S$, the union $R\cup S$ is defined to be the
graph, which contains crossings of both $R$ and $S$. If $R$ lies
above row $\ell$, while $S$ lies at or below row~$\ell$, then it
is easy to see $w(R\cup S)=w(R)w(S)$.

A "place $(i,j)$" of a graph $R$ will refer to either crossing or
non-crossing of strands in row $i$ and column $j$. For example, in
the second graph from Figure~\ref{more graphs} strands intersect
at place $(2,1)$, but do not intersect at place $(3,2)$. We refer
to those strands which intersect or do not intersect at place
$(i,j)$ as strands, which pass the place $(i,j)$. For instance,
strands $2,4$ pass places $(2,2)$ and $(1,3)$ in the third graph
of Figure~\ref{more graphs}. It will be convenient to write
$a\boxplus b=\ell$, if strand $a$ intersect strand $b$ in row
$\ell$ and strand $a$ is the horizontal strand of the crossing.
For example, $3\boxplus 4=3$ and $4\boxplus 3=1$ for the second
graph from Figure \ref{more graphs}.

\subsection{Schubert polynomials.}
For detailed discussions of Schubert polynomials $\mathfrak S_w$
we refer the reader to \cite{m-schubert} or \cite{man}. The only
property of Schubert polynomials used in this paper is stated in
Theorem \ref{def}, proved in \cite{f-s} and \cite{b-j-s}. So, for
purposes of this paper, we treat Theorem \ref{def} as a
definition.

For an rc-graph $R$ define $x^R=x_1^{|R_1|}x_2^{|R_2|}\dots$
(recall that in our notations $|R_i|$ is the number of crossings
of $R$ in row $i$).

\begin{theorem}
\label{def} For $w\in S_\infty$,
$$
\mathfrak S_w=\sum_{R\in\mathcal {RC}(w)}x^R.
$$
\end{theorem}

If $w$ is a shuffle $v(\lambda,r)$, then the Schubert polynomial
$\mathfrak S_w$ is known to be equal to the Schur polynomial
$S_\lambda(x_1,\dots, x_r)$ (for a definition of Schur polynomials
see~\cite{m}).

Schubert polynomials $\mathfrak S_w$ for all $w\in S_n$ form a
basis for $\mathbb C[x_1,\dots,x_n]$. Hence, for $u,v,w\in
S_\infty$, we can uniquely define the LR coefficients $c^u_{wv}$
by
$$
\mathfrak S_w\cdot\mathfrak S_v=\sum_{u}c^u_{wv} \mathfrak S_u.
$$

\subsection{Tableaux}
\begin{figure}[ht]
\begin{picture}(300,50)
\put (0,45){\line(1,0){15}} \put (0,30){\line(1,0){15}} \put
(15,45){\line(1,0){15}} \put (15,30){\line(1,0){15}} \put
(30,45){\line(1,0){15}} \put (30,30){\line(1,0){15}} \put
(00,15){\line(1,0){15}} \put (00,0){\line(1,0){15}} \put
(15,15){\line(1,0){15}} \put (15,15){\line(1,0){15}} \put
(30,15){\line(1,0){15}}

\put (0,30){\line(0,1){15}} \put (15,30){\line(0,1){15}} \put
(30,30){\line(0,1){15}} \put (45,30){\line(0,1){15}} \put
(00,15){\line(0,1){15}} \put (15,15){\line(0,1){15}} \put
(30,15){\line(0,1){15}} \put (45,15){\line(0,1){15}} \put
(00,0){\line(0,1){15}} \put (15,0){\line(0,1){15}}

\put(01,31) {\tiny 5}  \put (16,31) {\tiny 6} \put (31,31) {\tiny
7} \put (01,16) {\tiny 2} \put (16,16) {\tiny 3} \put (31,16)
{\tiny 4} \put (01,1) {\tiny 1}

\put (100,45){\line(1,0){15}} \put (100,30){\line(1,0){15}} \put
(115,45){\line(1,0){15}} \put (115,30){\line(1,0){15}} \put
(130,45){\line(1,0){15}} \put (130,30){\line(1,0){15}} \put
(100,15){\line(1,0){15}} \put (100,0){\line(1,0){15}} \put
(115,15){\line(1,0){15}} \put (115,15){\line(1,0){15}} \put
(130,15){\line(1,0){15}}

\put (100,30){\line(0,1){15}} \put (115,30){\line(0,1){15}} \put
(130,30){\line(0,1){15}} \put (145,30){\line(0,1){15}} \put
(100,15){\line(0,1){15}} \put (115,15){\line(0,1){15}} \put
(130,15){\line(0,1){15}} \put (145,15){\line(0,1){15}} \put
(100,0){\line(0,1){15}} \put (115,0){\line(0,1){15}}

\put(105,34) {1}  \put (120,34) {4} \put (135,34) {5} \put
(105,19) {3} \put (120,19) {5} \put (135,19) {6} \put (105,4) {4}

\put (200,45){\line(1,0){15}} \put (200,30){\line(1,0){15}} \put
(215,45){\line(1,0){15}} \put (215,30){\line(1,0){15}} \put
(230,45){\line(1,0){15}} \put (230,30){\line(1,0){15}} \put
(200,15){\line(1,0){15}} \put (200,0){\line(1,0){15}} \put
(215,15){\line(1,0){15}} \put (215,15){\line(1,0){15}} \put
(230,15){\line(1,0){15}}

\put (200,30){\line(0,1){15}} \put (215,30){\line(0,1){15}} \put
(230,30){\line(0,1){15}} \put (245,30){\line(0,1){15}} \put
(200,15){\line(0,1){15}} \put (215,15){\line(0,1){15}} \put
(230,15){\line(0,1){15}} \put (245,15){\line(0,1){15}} \put
(200,0){\line(0,1){15}} \put (215,0){\line(0,1){15}}

\put(203,34) {14}  \put (218,34) {26} \put (233,34) {27} \put
(203,19) {25} \put (218,19) {36} \put (233,19) {37} \put (203,4)
{35}

\end{picture}
\caption{Young diagram and tableaux} \label{young tableaux}
\end{figure}
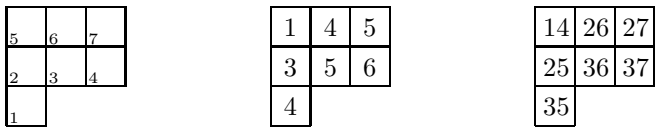

To a partition $\lambda=(\lambda_1\geq\dots\geq\lambda_r> 0)$
associate a Young diagram, which is given by $\lambda_i$ boxes in
row $i$. If $m=\sum^r_{i=1}\lambda_i$, label the boxes of the
Young diagram by integers $1$ to $m$ starting with the bottom row
going up and going from left to right in each row as shown in the
first picture of Figure \ref{young tableaux}, for
$\lambda=(3,3,1)$.

Given a Young diagram $D$, {\it a tableau} of shape $D$ is a
filling of boxes of the diagram $D$ by elements of certain
alphabets (we will talk later about the types of tableaux we
consider). For any tableau $U$, produce the word of $U$ by reading
the content of the boxes $1$ through $m$, denote it by $word(U)$.
For example, the word of the second tableau of Figure \ref{young
tableaux} is $4356145$. Set $|U|$ to be the number of boxes of
$U$.

Sometimes we consider only \emph{partially filled tableaux}: we
say that a tableau is \emph{filled up to $i$}, when boxes $1$ to
$i$ are filled, others are empty.

A tableau filled with positive integers is {\it row (column)
strict} if the numbers increase (do no decrease) from left to
right and do not decrease (increase) from top to bottom. For
example, second tableau of Figure \ref{young tableaux} is row and
column strict.

Given a shuffle $v(\lambda,r)$, by {\it a tableaux of
transpositions of $v(\lambda,r)$} we will understand a tableaux of
shape $\lambda$ filled by tuples $(ab)$ with $a\leq r<b$. For
example, see the third picture of Figure \ref{young tableaux},
where $\lambda=(3,3,1)$ and $r=3$. For a tableau of transpositions
$T$, let its word $word(T)$ be $(a_1b_1)\ldots(a_mb_m)$. Define
the permutations
$$
w_i(T)=t_{a_1b_1}\dots t_{a_{i}b_{i}}\text{ for }1\leq i\leq m,
\text{ and }w(T)=w_m(T).
$$
For $w\in S_\infty$, we say that $T$ is {\it an $r$-Brihat chain
of $w$}, if for $1\leq i\leq m$
\begin{equation}
\label{compatability} l(w w_i(T))=l(w)+i.
\end{equation}
If $T$ is filled up to $j$, we say that $T$ is an $r$-Bruhat chain
of $w$ if (\ref{compatability}) holds for $1\leq i\leq j$. For a
discussion of $r$-Bruhat chains see \cite{b-s}, where any sequence
$(a_1b_1)\ldots(a_mb_m)$, satisfying $l(wt_{a_1b_1}\dots
t_{a_{i}b_{i}})=l(w)+i$ and $a_i\leq r<b_i$ for all $i$, is an
$r$-Bruhat chain. For us it is convenient to think about
$(a_ib_i)$ as entries of a tableau.

We say that a triple $(w,R,T)$, consisting of a permutation $w$,
an rc-graph $R$ and a tableau of transpositions $T$ is an {\it
$r$-Bruhat package} if $ww(T)=w(R)$ and $T$ is an $r$-Bruhat chain
of $w$.

Associate to each permutation $w$ and each tableau of
transpositions $T$ of $v(\lambda,r)$ another tableau $E(w,T)$ of
the same shape. Fill box $i$ of $E(w,T)$ with $ww_i(T)(b_i)$. (If
$T$ is filled up to $j$, then $E(w,T)$ will be filled up to $j$.)
For example, if $w=(1,3,4,2,5,6,\dots)$ and $T$ is the third
tableaux from Figure \ref{young tableaux}, then the second tableau
from Figure \ref{young tableaux} is $E(w,T)$.

\subsection{Main Results}
\label{main-results} We are now ready to state our main results.
Given $R\in \mathcal {RC}(w)$ and $Y\in
\mathcal{RC}(v(\lambda,r))$ satisfying certain conditions, in
Section \ref{algorithm} we will define a graph $R\leftarrow Y$
together with a tableaux of transpositions $T(R,Y)$ of
$v(\lambda,r)$.

\begin{theorem}
\label{thm:main} Let $w,v(\lambda,r)\in S_{\infty}$ satisfy one
the following
\begin{align}
\label{case2} &\text{$w$ is an $r$-semi-shuffle,\ \ \ \ \ \ \ \ \
\ \ \ \ \ \ \ \ \ \ \ \ \ \ \ \ \ \ \ \ \ \ \ \ \ \ \ \ \ \ \ \ \ \ \ \ \ \ \ \ \ \ \ \ \ \ \ \ \ \ \ \ \ \ \ }\\
\label{case1} &\text{$\lambda$ is a hook.}
\end{align}
Let $R\in\mathcal{RC}(w)$ and $Y\in \mathcal{RC}(v(\lambda,r))$.
Then $U=R\leftarrow Y$ is an rc-graph and
\begin{align}
\label{thm:main-cond}&\text{$E(w,T(R,Y))$ is a
row and column strict tableau,\ \ \ \ \ \ \ \ \ \ \ \ \ \ \ \ \ \ \ \ \ \ \ \ \ \ \ \ \ \ \ }\\
\label{thm:cond1}&\text{$(w,U,T(R,Y))$ is an $r$-Bruhat package,}\\
\label{thm:cond3}&\text{$x^U=x^Rx^Y$.}
\end{align}
\end{theorem}

In Section \ref{inverse algorithm} the inverse insertion algorithm
will define graphs $U\rightarrow T$ and $Y(U,T)$ for certain
rc-graphs $U$ and tableaux of transpositions $T$.

\begin{theorem}
\label{thm:inverse main} Let $u,v(\lambda,r),w\in S_{\infty}$,
$U\in \mathcal{RC}(u)$ and $T$ be a tableau of transposition of
$v(\lambda,r)$. Assume $w,v(\lambda,r)$ satisfy (\ref{case1}) or
(\ref{case2}), $(w, U,T)$ is an $r$-Bruhat package and $E(w,T)$ is
row and column strict. Then $R=U\rightarrow T\in\mathcal{RC} (w)$
and $Y=Y(U,T)\in\mathcal {RC}(v(\lambda,r))$. Moreover,
$U=R\leftarrow T$ and $T=T(Y,R)$.
\end{theorem}

The next theorem is an immediate corollary of Theorems \ref{def},
\ref{thm:main}, \ref{thm:inverse main}.

\begin{theorem}
\label{thm:bruhat} Assume $w,u,v(\lambda,r)\in S_\infty$ satisfy
(\ref{case2}) or (\ref{case1}). Then $c^u_{wv(\lambda,r)}$ is
equal to the number of tableaux of transpositions $T$ of
$v(\lambda,r)$, such that $T$ is an $r$-Bruhat chain of $w$,
$E(w,T)$ is row and column strict and $ww(T)=u$.
\end{theorem}

Let us restate Theorem \ref{thm:bruhat} in the case when the shape
of $v$ is a hook in the form it appeared in \cite{s}. Given $w$
and a tableau of transpositions $T$ of $v$, define
$w^{(i)}=ww_i(T)$. If the shape of $v$ is a hook $(p, 1^{q-1})$
then $E(w,T)$ is row and column strict if and only if
\begin{equation}
\label{cond} w^{(1)}(b_1)>\dots >w^{(p)}(b_p) \text{ and
}w^{(p)}(b_p)<\dots <w^{(m)}(b_{m}),
\end{equation}
where $m=l(v)=p+q-1$. Using the fact that $t_{ab}$ and $t_{a'b'}$
commute as long as $a,b,a',b'$ are distinct, it can be shown that
there is a one to one correspondence between $r$-Bruhat chains,
which satisfy (\ref{cond}) and $r$-Brihat chains, which satisfy
\begin{equation}
\label{cond-a} w^{(1)}(a_1)>\dots >w^{(p)}(a_p) \text{ and
}w^{(p)}(a_p)<\dots <w^{(m)}(a_{m}).
\end{equation}
(This correspondence can be constructed by starting with a chain,
which satisfies (\ref{cond}) and commuting transpositions of this
chain to make sure (\ref{cond-a}) holds.)

So Theorem \ref{thm:bruhat} in case (\ref{case1}) can be restated
as it originally appeared in \cite{s}.

\begin{theorem}\label{thm:hook} Assume $v=((p,1^{q-1}),r)$. Then $\mathfrak
S_w\mathfrak S_v=\sum \mathfrak S_{w^{(m)}}$, the sum over all
paths in $k$-Bruhat order $w<_rw^{(1)}<_r\dots<_rw^{(m)}$, which
satisfy \ref{cond-a}.
\end{theorem}

\section{Insertion Algorithm}
\label{algorithm}

\subsection{Preliminaries} We need some preliminaries before defining the algorithm.

First, let $Y$ be an rc-graph with $w(Y)=v(\lambda,r)$. Sometimes
we will think of the Young diagram of $\lambda$ as the shape of
$Y$, denoted by $sh(Y)$. It is easy to see that strand $s\leq r$
intersects exactly $\lambda_{r+1-s}$ other strands. Let these
intersections be in rows $i_1\geq\dots\geq i_{\lambda_{r+1-s}}$
(one number for each crossing, so that repetitions are allowed),
then define $word(Y,s)=i_1\dots i_{\lambda_{r+1-s}}$. Define
$word(Y)$ to be the concatenation $word(Y,1)\dots word(Y,r)$.
Notice that if two strands $a,b$ intersect in $Y$ and $a<b$ then
$a\leq r<b$. Hence every crossing of $Y$ correspond to a single
letter in $word(Y)$. Also notice that for any $\ell$, the
permutation of $Y_{\geq \ell}$ is again a shuffle. Moreover, the
shape of $Y_{\geq \ell}$ is a subdiagram of the shape of $Y$.

Secondly, we need the following lemma and the construction after
the lemma.

\begin{lemma}
\label{add-crossing} (1) If $R$ is an rc-graph and
$l(w(R)t_{cd})=l(w(R))-1$, then strands $c$ and $d$ intersect in
$R$, and removing this crossing produces another rc-graph.

\noindent (2) Let $R$ be an rc-graph with strands $c$ and $d$
passing place~$(\ell,j)$ but never intersecting in $R$. Then
insertion a crossing into place $(\ell,j)$ produces an rc-graph.
\end{lemma}

\proof For $w\in S_\infty$ its length is $l(w)=\#\{(i,j):i<j,
w(i)>w(j)\}$. Hence
\begin{align}
\label{no i}
&\text{for $c<d$, $l(wt_{cd})=l(w)+1$ if and only if $w(c)<w(d)$ and}\\
\nonumber &\text{there is no $i$ with } c<i<d, \ \ w(c)<w(i)<w(d).
\end{align}

To prove the first part of the lemma, notice that (\ref{no i})
applied to $w=w(R)t_{cd}$ immediately implies that
$w(R)(c)>w(R)(d)$. In particular, strands $c$ and $d$ must
intersect in $R$. Remove their crossing to produce graph $R'$.
Using (\ref{number-length}), (\ref{simple-add-crossing}) and
$l(w(R'))=l(w(R)t_{cd})=l(w(R))-1=|R'|$, we conclude $R'$ is an
rc-graph.

To prove the second part, add the crossing of strand $c$ and $d$
in place $(\ell,j)$ of $R$ to produce graph $R'$. By
(\ref{number-length}) it is enough to check $l(w(R'))=l(w(R))+1$.
Since strands $c$ and $d$ do not intersect in $R$,
$w(R)(c)<w(R)(d)$, hence by (\ref{no i}), it is enough to check
that there is no $i$ with
$$c<i<d, \,\,\,
w(R)(c)<w(R)(i)<w(R)(d).
$$
It is very easy to see that if such $i$ existed, strand $i$ would
have to intersect either strand $c$ or strand $d$ twice in $R$,
which is impossible.
\endproof

As a consequence to Lemma \ref{add-crossing} let us present the
following construction. Given an $r$-Bruhat package $\mathcal
P=(w,R,T)$,  let $m=|T|$. Set $S_m(\mathcal P)=R$. Then, by Lemma
\ref{add-crossing}, rc-graphs $S_j(\mathcal P)$ for $0\leq j\leq
m$ are uniquely defined, once we require that $w(S_j(\mathcal
P))=ww_j(T)$, and $S_{j}(\mathcal P)$ is constructed out of
$S_{j+1}(\mathcal P)$ by removing exactly one crossing. An example
is provided in Section~\ref{example1}.

\subsection{Outline}
We now start the description of the algorithm. We begin with a
general outline to better explain the procedure.

Assume we are given an rc-graph $R$ with $w=w(R)$ and an rc-graph
$Y$ with $v(\lambda,r)=w(Y)$, satisfying (\ref{case2}) or
(\ref{case1}). Our goal is to define a graph $R\leftarrow Y$ and a
tableau of transpositions $T(R,Y)$ of $v(\lambda,r)$.

The algorithm starts with row $r$ and goes up. After the algorithm
finishes with row $\ell$ it produces
$R(\ell)=R_{\geq\ell}\leftarrow Y_{\geq\ell}$  and
$T(\ell)=T(R_{\geq\ell},Y_{\geq\ell})$.

It is convenient to think of $T(\ell)$ as "the history" of the
algorithm up to $\ell$. Namely, $T(\ell)$ says how to go from
$w(R_{\geq\ell})$ to $w(R(\ell))$ along a chain in Bruhat order.
Moreover, since the shape of $Y_{\geq\ell}$ is the same as the
shape of $T(\ell)$, $word(Y_{\geq\ell})$ and $word(T(\ell))$ are
of the same shape, so each step in the chain corresponds to a
crossing of $Y_{\geq\ell}$.

Assume the insertion has been performed up to row $\ell+1$. The
next row where it has to operate is row $\ell$. Clearly
$word(Y_{\geq\ell})$ is constructed from $word(Y_{\geq\ell+1})$ by
adding letter $\ell$ at some places. If
$word(Y_{\geq\ell})=word(Y_{\geq\ell+1})$, so there are no
crossings in row $\ell$ of $Y$, then it looks like we can say
$T(\ell)=T(\ell+1)$. But since $T(\ell+1)$ defines a chain which
starts at $w(R_{\geq\ell+1})$, it may not define a chain starting
from $w(R_{\geq \ell})$. So, the algorithm goes through the
rectification of both rc-graph and the chain to make sure
$T(\ell)$ indeed defines the chain starting at $w(R_{\geq\ell})$.

If there are some crossings in row $\ell$, then the algorithm
inserts them whenever necessary. The order in which insertions and
rectifications are performed are determined by the order of letter
in $word(Y_{\geq\ell})$.

\subsection{Sequence of Steps of the algorithm.}
\label{steps} Throughout the rest of this section all statements,
which require proofs, will be underlined and then proved in
Section~\ref{proof}.

For each $(\ell, i)$ with $r\geq \ell \geq 1$ and $0\leq i\leq
m_\ell$ (where $m_\ell=|Y_{\geq \ell}|$) the algorithm performs a
step, which we call \emph{step $(\ell,i)$}. The steps go in the
following order. Step $(\ell,i+1)$ goes after step $(\ell,i)$ if
$m_\ell-1\geq i\geq0$. Step $(\ell,0)$ follows step
$(\ell+1,m_{\ell+1})$.

Before giving a detailed description of each step, let us present
the data produced by each step and the conditions this data
satisfies. Step $(\ell,i)$ constructs
\begin{align}
&\text{rc-graphs $R(\ell,i)$,
with no crossings above row $\ell$,} \\
&\text{tableau of transpositions $T(\ell,i)$ for the shuffle
$w(Y_{\geq \ell})$ filled up to $i$. \ \ \ \ \ \ \ \ }
\end{align}
Here, $R(\ell,i)$ and $T(\ell,i)$ play the role of the
intermediate result of the algorithm. After each step, $R(\ell,i)$
and $T(\ell,i)$ must satisfy the following conditions
\begin{align}
\label{main-cond}&\text{$E(w(R_{\geq\ell}),T(\ell,i))$ is a
row and column strict tableau,}\\
\label{cond1}&\text{$\mathcal P(\ell,i)=(w(R_{\geq\ell}),
R(\ell,i), T(\ell,i))$ is an $r$-Bruhat package. \ \ \ \ \ \ \ \ \
\ \ \ \ \ \ \ \ \ \ \ \ \ }
\end{align}

\begin{remark}
{\rm  Conditions (\ref{main-cond}), (\ref{cond1}) are analogues of
(\ref{thm:main-cond}), (\ref{thm:cond1}) from
Theorem~\ref{thm:main}. That is why they have to be satisfied by
the intermediate results of the algorithm. Condition
(\ref{main-cond}) is the condition, which we do not know how to
generalize to cases other than (\ref{case2}) and (\ref{case1}).}
\end{remark}

So for each row $\ell$ we start with {\it a row-to-row} step
$(\ell,0)$, which sets up the data needed for performing insertion
in this row. Then we perform a step for each letter of
$word(Y_{\geq\ell})$. If this letter is equal to $\ell$ it {\it an
insertion} step and we will insert a crossing in row $\ell$ to the
current rc-graph. If the letter is not $\ell$, then we perform
{\it a rectification} and rectify, if necessary, both rc-graph
$R(\ell,i)$ and the chain given by $T(\ell,i)$ to guarantee both
(\ref{main-cond}) and (\ref{cond1}) are satisfied.

The rest of Section \ref{steps} introduces additional notation and
states two additional conditions, which clarify certain parts of
the algorithm and simplify certain proofs.

After we are finished with all the steps for row $\ell$, we are
given $R(\ell,m_\ell)$ and $T(\ell,m_\ell)$, which, to shorten the
notations, we denote by $R(\ell)$ and $T(\ell)$.  Denote by
$\mathcal P(\ell)$ the $r$-Bruhat package $(w(R_{\geq\ell}),
R(\ell), T(\ell))$.

Fix $\ell$, let $word(Y_{\geq \ell})=k_1\dots k_{m_\ell}$. Each
letter $k_i$ of $word(Y_{\geq \ell})$ corresponds to a crossing of
$Y_{\geq\ell}$ in row $k_i$. If $k_i>\ell$, then the letter $k_i$
is also a part of $word(Y_{\geq \ell+1}) =k'_1\dots
k'_{m_{\ell+1}}$, let the index of $k_i$ inside $word (Y_{\geq
\ell+1})$ be $i_+$. Set $i_+=0$, if $i=0$. So, if we think of
$sh(Y_{\geq\ell+1})$ as a subdiagram of $sh(Y_{\geq \ell})$, then
box~$i$ of $Y_{\geq \ell}$ coincides with box $i_+$ of
$sh(Y_{\geq\ell+1})$. The two additional conditions are
\begin{align}
\label{cond3}&\text{$x^{R(\ell)}=x^{R_{\geq\ell}}x^{Y_{\geq\ell}},$}\\
\label{cond4}&\text{If $k_i<\ell$, then $R(\ell,i)_{\geq \ell+1}
=S_{i_+}(\mathcal P(\ell+1))$.\ \ \ \ \ \ \ \ \ \ \ \ \ \ \ \ \ \
\ \ \ \ \ \ \ \ \ \ \ \ \ \ \ \ \ \ \ \ \ \ \ \ }
\end{align}
It will be obvious from the description of the algorithm that
these conditions are always satisfied. Condition (\ref{cond3})
implies that (\ref{thm:cond3}) holds for the final result, while
(\ref{cond4}) indicates, that step $(\ell,i)$ only operates in row
$\ell$, as the part of $R(\ell,i)$, which lies below row $\ell$,
is uniquely determined by $\mathcal P(\ell+1)$.

\subsection{Start of the algorithm} Set $R(r,0)=R_{\geq r}$ and
let $T(r,0)$ be the empty tableau of shape $sh(Y_{\geq r})$, then
\underline{$R(r,0)$ and $T(r,0)$ satisfy (\ref{main-cond}),
(\ref{cond1})}.

\subsection{Row-to-row steps} Each step $(\ell,0)$ is called
\emph{a row-to-row step}. This step sets $T(\ell,0)$ to be the
empty tableau of shape $sh(Y_{\geq\ell})$ and
\begin{equation}
\label{eq:r(ell,0)} R(\ell,0)=S_0(\mathcal P(\ell+1))\cup R_\ell.
\end{equation}
Then \underline{$R(\ell,0)$ is an rc-graphs} and
\underline{$R(\ell,0)$ and $T(\ell,0)$ satisfy (\ref{main-cond}),
(\ref{cond1})}.

As mentioned before, this step sets up the data for performing the
algorithm in row $\ell$. $T(\ell+1)$ defines a chain from
$w(R_{\geq\ell+1})$ to $w(R(\ell+1))$. On the level of rc-graph
this chain is given by the chain $S_0(\mathcal P(\ell+1)),\dots,
S_{m_{\ell+1}}(\mathcal P(\ell+1))$. So we can thing of
$R(\ell,0)$ as backtracking the algorithm from $R(\ell+1)$ to
$S_0(\mathcal P(\ell+1))$ and then adding those crossings of $R$,
which lie in row $\ell$.

\subsection{Insertions}
Assume $word(Y_{\geq\ell})=k_1\dots k_{m_\ell}$. If $k_i$ is the
letter $\ell$, then step $(\ell,i)$ is called \emph{an insertion}.

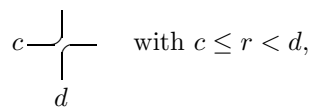
\begin{figure}[ht]
\begin{picture}(100, 25)

\put (0,15){\line(1,0){10}} \put (13,2){\line(0,1){10}} \put
(16,15){\line(1,0){10}} \put (13,18){\line(0,1){10}}

\put(9,19){\oval(8,8)[br]} \put(17,11){\oval(8,8)[tl]}

\put(-6,13){$c$} \put(10,-8){$d$}

\put(40,13){with $c\leq r <d$,}

\end{picture}
\caption{Place where insertion is allowed.} \label{allowed}
\end{figure}

During insertion step $(\ell,i)$, we say that insertion into a
place $(\ell,j)$ is allowed, if strands $c,d$ pass this place in
$R(\ell,i-1)$ as shown in Figure \ref{allowed} and $c\leq r<d$.

Find the rightmost place, where insertion is allowed. Let it be
place $(\ell,j_0)$ and let strands $c$ and $d$ pass through it.
Define $R(\ell,i)$ by adding a crossing to $R(\ell,i-1)$ into
place $(\ell,j_0)$. Define $T(\ell,i)$ by adding $(cd)$ to box $i$
of $T(\ell,i-1)$. Then \underline{$R(\ell,i)$ is an rc-graph} and
\underline{(\ref{main-cond}) and (\ref{cond1}) are satisfied.}

\subsection{Rectifications }
If $k_i>\ell$, then step $(\ell,i)$ is called {\it a
rectification}. The first part of rectification is to define a
graph $R'$ and a tableau of transpositions $T'$. The rc-graph
$S_{i_+}(\mathcal P(\ell+1))$ has one more crossing than
$S_{i_+-1}(\mathcal P(\ell+1))$, add this crossing to
$R(\ell,i-1)$ to produce $R'$. Then, since (\ref{cond4}) holds for
$R(\ell,i-1)$, $R'$ coincides with $S_{i_+}(\mathcal P(\ell+1))$
below row $\ell$ and row $\ell$ of $R'$ is the same as row $\ell$
of $R(\ell,i-1)$. To produce $T'$, add to box $i$ of $T(\ell,i-1)$
the entry $(ab)$ of box $i_+$ of $T(\ell+1)$.

\begin{figure}[ht]
\begin{picture}(330, 25)

\put (0,15){\line(1,0){28}} \put (13,2){\line(0,1){28}}

\put(-6,13){$b$} \put(10,-8){$a$}

\put(40,20){where $(ab)$ is the} \put(40,5){entry of box $i$ of
$T'$.}


\put (150,15){\line(1,0){28}} \put (163,2){\line(0,1){28}}

\put(140,13){$b$} \put(160,-8){$f$}

\put(190,20){where $(ab)$ and $(af)$ are the} \put(190,5){entries
of boxes $i$ and $i-1$ of $T'$.}

\end{picture}

\caption{One of these crossing in row $\ell$ needs to be removed.}
\label{remove-crossing}
\end{figure}
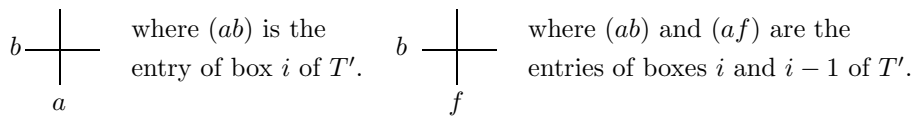

If $E(w(R_{\geq\ell}), T')$ is row and column strict and
$(w(R_{\geq\ell}), R',T')$ is an $r$-Bruhat package, set
$R(\ell,i)=R'$ and $T(\ell,i)=T'$ and move on to the next step.
Otherwise, \underline{there is a crossing in $R'$ in row $\ell$,
which fits the description in Figure \ref{remove-crossing}.}

If $R'$ has a crossing, which looks like the first crossing of
Figure \ref{remove-crossing}, remove this crossing to produce
$R''$ and remove $(ab)$ from box $i$ of $T'$ to produce $T''$.
Otherwise remove the crossing of strands $b$ and $f$ shown in
Figure \ref{remove-crossing} to produce $R''$, remove $(ab)$ from
box $i$ in $T'$ and replace the entry of box $i-1$ of $T'$ by
$(ab)$ to produce~$T''$.

\begin{figure}[ht]
\begin{picture}(320, 25)

\put (0,15){\line(1,0){10}} \put (13,2){\line(0,1){10}} \put
(16,15){\line(1,0){10}} \put (13,18){\line(0,1){10}}

\put(9,19){\oval(8,8)[br]} \put(17,11){\oval(8,8)[tl]}

\put(-6,13){$c$} \put(10,-8){$d$}

\put(40,13){with $c\leq r <d$,}


\put (150,15){\line(1,0){10}} \put (163,2){\line(0,1){10}} \put
(166,15){\line(1,0){10}} \put (163,18){\line(0,1){10}}

\put(159,19){\oval(8,8)[br]} \put(167,11){\oval(8,8)[tl]}

\put(140,13){$g$} \put(160,-8){$d$}

\put(190,20){where $(eg)$ is the entry} \put(190,5){of box $i-1$
of $T''$, and $r<d$}

\end{picture}

\caption{Places, where insertions are allowed during
rectification.} \label{allowed-insertions-rect}
\end{figure}
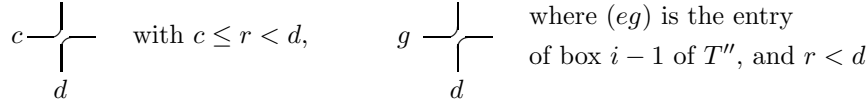

We say that insertions into places in $R''$ shown in
Figure~\ref{allowed-insertions-rect} are allowed.

Remember that $R''$ is constructed out of $R'$ by removing a
crossing from place $(\ell,j_0)$. Find rightmost place to the left
of $(\ell,j_0)$, where insertion is allowed. Insert a crossing
there to produce $R(\ell,i)$. If this is the place of the first
type from Figure \ref{allowed-insertions-rect}, then insert $(cd)$
into box $i$ of $T''$ to define $T(\ell,i)$. In the second case,
replace the entry of box $i-1$ of $T''$ by $(ed)$ and place $(eg)$
into box $i$ of $T''$ to define $T(\ell,i-1)$. Then
\underline{$R(\ell,i)$ is an rc-graph} and
\underline{(\ref{main-cond}) and (\ref{cond1}) are satisfied.}

\subsection{End of the algorithm.} Set $R\leftarrow Y=R(1,m_1)$
and $T(R,Y)=T(1,m_1)$.

\subsection{Concluding remarks} As mentioned before, the algorithm
can be simplified in all the cases, except for the case when shape
of $v$ is a hook. For case (\ref{case2}) (see \cite{k}) our
algorithm produces the same result as inserting letters of
$word(Y)$ one by one into $R$ using algorithm of Bergeron and
Billey \cite{b-b}. Moreover, in the case when $w$ is also an
$r$-shuffle, it is just the Schensted insertion algorithm. So our
algorithm is a direct generalization of Schensted insertion
algorithm for Young tableaux.

For the case when shape of $v$ is a row a simplified algorithm is
given in \cite{k-k}. In the case when shape of $v$ is a column,
the only simplification, which can be done, is omitting the second
picture from Figure \ref{allowed-insertions-rect}. Analogous
algorithm in this case can also be found in \cite{ku}.

\section{Proof of Theorem \ref{thm:main}.}
\label{proof} To prove Theorem \ref{thm:main}, it is enough to
prove all the statements underlined in Section \ref{algorithm}.
Let us list these statements again.
\begin{enumerate}
\item{$R(r,0)$ and $T(r,0)$ satisfy conditions
(\ref{main-cond}), (\ref{cond1}).}
\item{For $r>\ell\geq 1$, $R(\ell,0)$ is an rc-graph and $R(\ell,0)$, $T(\ell,0)$ satisfy
(\ref{main-cond}), (\ref{cond1}).}
\item{After insertion, $R(\ell,i)$ is an rc-graph and (\ref{main-cond}), (\ref{cond1}) are satisfied.}
\item{During rectification, if $R'$ is not an rc-graph or $E(w(R_{\geq\ell},T')$ is not row and column strict,
there is a crossing in $R'$ shown in Figure \ref{remove-crossing}.}
\item{After rectification, $R(\ell,i)$ is an rc-graph and (\ref{main-cond}), (\ref{cond1}) are satisfied.}
\end{enumerate}

\subsection{Proof of (1)} Since $T(r,0)$ is an empty tableau,
condition (\ref{main-cond}) is vacuous, while (\ref{cond1})
follows directly from $R(r,0)=R_{\geq r}$.

\subsection{Proof of (2)} Since $S_0(\mathcal P(\ell+1))$  has no
crossings above row $\ell+1$ and $R_\ell$ has only crossings in
row $\ell$, we know from (\ref{eq:r(ell,0)})
$$
w(R(\ell,0))=w(S_0(\mathcal P(\ell+1)))w(R_\ell)
=w(R_{\geq\ell+1})w(R_\ell) =w(R_{\geq\ell+1}\cup R_\ell)
=w(R_{\geq\ell}).
$$
On the other hand
$$
|R(\ell,0)|= |S_0(\mathcal P(\ell+1))| +|R_\ell|
=|R_{\geq\ell+1}|+|R_\ell|=|R_{\geq\ell}|.
$$
Hence $l(w(R(\ell,0))) =|R(\ell,0)|$ and using
(\ref{number-length}) we conclude $R(\ell,i)$ is an rc-graph.

Since $T(r,0)$ is an empty tableau, condition (\ref{main-cond}) is
vacuous, while (\ref{cond1}) follows immediately from
$w(R(\ell,0))=w(R_{\geq\ell})$.

\subsection{Proof of (3)}
To show that the algorithm works properly during insertion, we
must show that there are places in $R(\ell,i-1)$, where insertions
are allowed. To do this, define the sequence $c_k$ by $c_0=\ell$
and
\begin{equation}
\label{sequence ck} a_k=w(R(\ell,i-1))(c_k), \,\,\,
c_{k+1}=w(R(\ell,i-1)_{\geq\ell+1})^{-1}(a_{k}+1).
\end{equation}
Here is another way of defining $c_k$'s. Look at all strands of
$R(\ell,i-1)$, which have horizontal parts in row~$\ell$, in other
words, which do not cross row $\ell$ vertically. These strands do
not cross each other in row $\ell$. So, we let the sequence $c_k$
be the labels of these strands, read from left to right. For
example, if $R(\ell,i-1)$ is given by the third graph of
Figure~\ref{more graphs} and $\ell=1$, then $c_0=1$, $c_1=4$,
$c_2=2$ and so on.

It is clear that for large $k$, $c_k>r$. Since $c_0=\ell\leq r$,
there exists $\bar k$ with $c_{\bar k}\leq r< c_{\bar k+1}$. By
construction, we know that strands $c_{\bar k}$ and $c_{\bar k+1}$
pass next to each other in row $\ell$ in some place $(\ell,j)$.
Then insertion into place $(\ell,j)$ is allowed, as $c=c_{\bar
k}\leq r< c_{\bar k}=d$ as required in Figure \ref{allowed}.

\subsubsection{$R(\ell,i)$ is an rc-graph and (\ref{cond1}) holds.}
$R(\ell,i)$ is an rc-graph by the second part of Lemma
\ref{add-crossing}.  Moreover,
$l(w(R(\ell,i))=l(w(R(\ell,i-1)))+1$ and
\begin{equation}
\label{c-d} w(R(\ell,i))= w(R(\ell,i-1))t_{cd}, \ \
w(T(\ell,i))=w(T(\ell,i-1))t_{cd},
\end{equation}
which immediately implies that (\ref{cond1}) holds.

\subsubsection{After insertion, (\ref{main-cond}) holds.}
If $i=1$, so that the insertion step corresponds to the first
letter $k_1=\ell$, condition (\ref{main-cond}) is vacuous.
Otherwise, we will show that there exist $j$, such that insertion
into $(\ell,j)$ is allowed and
\begin{equation}
\label{check2} w(R(\ell,i-1)t_{c'd'})(d')>w(R(\ell,i-1))(f)
\end{equation}
where $c',d'$ are the strands passing place $(\ell,j)$ and $(ef)$
is the entry of box $i-1$ of $T(\ell,i-1)$. This will be enough to
prove (\ref{main-cond}). Indeed if $R(\ell,i)$ is defined by
adding a crossing of strand $c,d$ in place $(\ell,j_0)$, then
$j_0\geq j$ and therefore (\ref{check2}) holds for  $c',d'$
substituted by $c,d$. Hence $E(w(R_{\geq\ell}),T(\ell,i))$ is row
and column strict.

To show that $j$, satisfying (\ref{check2}), exists, consider
sequence $c_k$ with different $c_0$: if strand $e$ in
$R(\ell,i-1)$ intersects vertically another strand $e'$ in row
$\ell$, that is $e'\boxplus e=\ell$, then set $c_0=e'$ otherwise
set $c_0=e$ (notice that $e'<e\leq r$). Since $c_0\leq r$ and
$c_k>r$ for large $k$, there exists $\bar k$ with $c_{\bar k}\leq
r< c_{\bar k+1}$. Let $c'=c_{\bar k}$ and $d'=c_{\bar k+1}$. Then
strands $c'$ and $d'$ pass next to each other in row $\ell$ at a
place $(\ell,j)$ and insertion into $(\ell,j)$ is allowed.
Moreover, since strand $c_0$ is either strand $e$ or it intersects
strand $e$ horizontally in row $\ell$, the following calculation
proves (\ref{check2})
$$
w(R(\ell-1,i)t_{c'd'})(d')=w(R(\ell,i-1))(c')\geq
w(R(\ell,i-1))(e)>w(R(\ell,i-1))(f).
$$

\subsection{Proof of (4)} Recall that $R'$ is constructed out of
$R(\ell,i-1)$ by adding a crossing of strands $a,b$ to guarantee
$R'$ coincides with $S_{i_+}(\mathcal P(\ell+1))$ below row
$\ell$.

Let us show that $R'$ is not an rc-graph if and only if strands
$a$ and $b$ intersect in row $\ell$ as shown in Figure
\ref{remove-crossing}. Indeed, if $a,b$ intersect in row $\ell$,
they intersect twice in $R'$, so $R'$ is not an rc-graph.
Conversely, if they do not intersect in row $\ell$ of $R'$, then
they do not intersect in $R(\ell,i-1)$ (If they intersect below
row $\ell$ then $S_{i_+}(\mathcal P(\ell+1))$ is not an rc-graph.)
Hence by the second part of Lemma~\ref{add-crossing}, $R'$ is an
rc-graph.

Since $w(R')=w(R(\ell,i-1))t_{ab}$, we can immediately conclude
that if $R'$ is an rc-graph, then $(w(R_{\geq\ell}),R', T')$ is an
$r$-Bruhat package.

To finish the proof of (4), it remains to prove that if
$E(w(R_{\geq\ell}), T')$ is not column or row strict then the
second crossing from Figure \ref{remove-crossing} must occur. We
will do it separately for cases (\ref{case2}) and (\ref{case1}).

\subsubsection{Case (\ref{case2})}
\label{proof(4),case1} We start with preliminaries, which we also
use in the proof of~(5).

\begin{lemma}
\label{case2_holds} If $u\in S_\infty$ is an $r$-semi-shuffle,
$u'=ut_{ab}$ with $a\leq r<b$ and $l(u')=l(u)+1$, then $u'$ is an
$r$-semi-shuffle .
\end{lemma}

\proof We must show that if $r<b'<b''$ then
$ut_{ab}(b')<ut_{ab}(b'')$. If $b'\neq b$ and $b''\neq b$, then
$ut_{ab}(b')=u(b')<u(b'')=ut_{ab}(b''),$ since $u$ is an
$r$-semi-shuffle.

If $b'=b$, then $ut_{ab}(b')=u(a)<u(b)<u(b'')=ut_{ab}(b''),$ since
$l(ut_{ab})=l(u)+1 $ and $u$ is an $r$-semi-shuffle.

If $b''=b$, then $ut_{ab}(b')=u(b')<u(a)=ut_{ab}(b''),$ where
$u(b')<u(a)$, since otherwise $a<b'<b$ and $u(a)<u(b')<u(b)$,
which contradicts (\ref{no i}).
\endproof

Assume  $T$ is a tableau of transpositions (possibly partially
filled). If $(a_kb_k)$ are the entries of $T$, let $B(T)$ be the
tableau of the same shape with the entries $b_k$.

\begin{lemma}
\label{lemma:standard} Let $w$ be an $r$-semi-shuffle and $T$ be
an $r$-Bruhat chain of $w$. Then $E(w,T)$ is row and column strict
if and only if $B(T)$ is row strict.
\end{lemma}

\proof Assume $u\in S_{\infty}$ is an $r$-semi-shuffle and
$l(ut_{ab})=l(u)+1$ for $a\leq r<b$. Let $b'>r$ then it is easy to
see by Lemma \ref{case2_holds} that
\begin{equation}
\label{eq:check5} \text{$u(b')<ut_{ab}(b)$ if and only if $b'<b$.}
\end{equation}
Clearly, (\ref{eq:check5}) implies that rows of $E(w,T)$ strictly
increase from left to right if and only if the same holds for
$B(T)$.

Let us show that if $B(T)$ is row strict, then $E(w,T)$ is column
strict. Denote by $e_k$ the entry of box $k$ of $E(w,T)$. Let box
$i'$ be directly above box $i$ in $sh(T)$. Consider boxes $i$
through $i'$ in the diagram $sh(T)$, as shown in Figure \ref{row
of boxes}.

\begin{figure}[ht]
\begin{picture}(200,60)
\put(0,0){\line(1,0){96}} \put(0,0){\line(0,1){60}}
\put(96,12){\line(1,0){36}} \put(96,0){\line(0,1){12}}
\put(132,48){\line(1,0){24}} \put(132,12){\line(0,1){36}}
\put(0,60){\line(1,0){156}} \put(156,48){\line(0,1){12}}

\put(0,36){\line(1,0){132}} \put(0,48){\line(1,0){60}}
\put(48,24){\line(1,0){84}}

\put(12,36){\line(0,1){12}} \put(48,24){\line(0,1){24}}
\put(60,24){\line(0,1){24}} \put(120,24){\line(0,1){12}}

\put(01,38) {\tiny $j$} \put(49,37) {\tiny $i'$} \put(49,25)
{\tiny $i$} \put(120,26) {\tiny$j$\tiny$-$\tiny$1$}
\end{picture}
\caption{Boxes $i$ through $i'$ of $T$.} \label{row of boxes}
\end{figure}
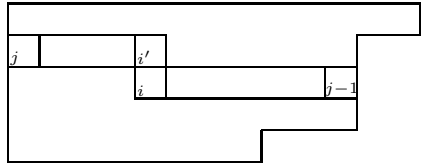

To show $E(w,T)$ is column strict, it is enough to show
$e_{i'}<e_i$, for any $i$ not in the top row. If $B(T)$ is
row-strict, then $b_{i}\neq b_{\tilde i}$ for $i<\tilde i<i'$. If
$b_{i'}<b_i$, then
$$
e_i=ww_i(T)(b_i)= ww_{i'}(T)(b_i)>ww_{i'}(T)(b_{i'})=e_{i'},
$$
since by Lemma \ref{case2_holds} $ww_{i'}(T)$ is an
$r$-semi-shuffle. If $b_{i'}=b_{i}$, then
$$
e_i=ww_i(T)(b_i)= ww_{i'-1}(T)(b_i)=
ww_{i'}(T)(a_i)>ww_{i'}(T)(b_i)=e_{i'},
$$
since $l(ww_{i'-1}(T))+1=l(ww_{i'}(T))$.

Conversely, assume $E(w,T)$ is row and column strict. To show
$B(T)$ is row strict, it is enough to show $b_{i'}\leq b_i$ for
any $i$, such that box $i$ is not in the top row. Assume for a
moment $b_{i}\neq b_{\tilde i}$ for any $j\leq \tilde i\leq i'$.
Then
$$ww_{i'}(T)(b_{i})=ww_{i}(T)(b_i)=e_i>e_{i'}= ww_{i'}(T)(b_{i'}),$$
since $E(w,T)$ is row and column strict. Thus, since $ww_{i'}(T)$
is an $r$-semi-shuffle, we conclude $b_{i'}<b_i$.

Otherwise, if $b_i=b_{\tilde i}$ for some $j\leq \tilde i\leq i'$,
to show that $\tilde i=i'$, use induction on $i$. If $i$ is the
first box in its row, then $b_{i'}= b_i$, as $i'=j$. Otherwise,
assume the box underneath box $\tilde i$ contain $\bar b$. By
induction $\bar b\geq b_{\tilde i}$. On the other hand, we know
that $\bar b<b_{i}$, if $\tilde i_{-}\neq i$. Hence, if
$b_i=b_{\tilde i}$, then box $i$ must be underneath box $\tilde
i$.
\endproof

To finish the proof of (4) in case (\ref{case2}), we will prove
that if $R'$ is an rc-graph, then $E(w(R_{\geq\ell}),T')$ is row
and column strict, or by Lemma \ref{lemma:standard} it is enough
to show $B(T')$ is row strict.

Let $b_k$ be the entries of $B(T')$, so that $b_i=b$. Since boxes
$1$ through $i-1$ of $B(T')$ and $B(T(\ell,i-1))$ coincide,
$B(T')$ can fail to be row strict if box $i-1$ is in the same row
as box $i$ and $b_{i-1}\geq b_i=b$, or, if there is box $i_*$
underneath box $i$, such that $b=b_i> b_{i_*}$. Let us show both
cases are impossible. This will finish the proof of (4) in case
(\ref{case2}).

If $(\bar e\bar f)$ is the entry of box $i_+-1$ of $T(\ell+1)$,
then by construction $b_{i-1}\leq \bar f$. At the same time if box
$i-1$ is in the same row as box $i$, then $\bar f<b_i=b$, since
$B(T(\ell+1))$ is row strict. So $b_{i-1}<b$, whenever box $i-1$
is in the same row as box $i$.

It remains to show that if box $i_*$ is the box underneath box $i$
in $B(T')$ , then $b\leq b_{i^*}$. We will prove it by induction
on $i$. We will prove the induction step when step $(\ell+1,i_*)$
is an insertion. If it is a rectification, the proof is almost
identical.

Denote temporarily $\tilde R=R(\ell,i_*-1)$. We will prove there
exists a place in $\tilde R$ in row $\ell+1$ shown in Figure
\ref{allowed} with $d\geq b$. Then we will be guaranteed
$b_{i_*}\geq b$.

Look at how stand $b$ passes row $\ell+1$ in $\tilde R$. If it
passes it vertically, then it intersects certain strand $a'$ with
$a'\leq r$ (since $w(\tilde R)$ is an $r$-semi-shuffle). Then
consider the sequence $c_k$ for rc-graph $\tilde R$, as defined in
(\ref{sequence ck}), with $c_0=a'$. By the same argument as in the
proof of (3) we can find a place $(\ell,j)$, shown in Figure
\ref{allowed}, to the right of the place where strand $b$ passes
row $\ell+1$. Hence $d>b$.

If strand $b$ does not pass row $\ell$ vertically, look again at
the sequence $c_k$ for $\tilde R$ with $c_0=\ell$. Strand $b$ is
an element of this sequence, let $b=c_{\tilde k}$. Let us show
that
\begin{equation}
\label{eq:check3} c_{\tilde k-1}\leq r\ \ \ \text{ or }\ \ \
c_{\tilde k-1}=b-1.
\end{equation}
Indeed, if $c_{\tilde k-1}>r$, then if there exist $b'$ with
$c_{\tilde k-1}<b'<b$, then strand $b'$ must intersect either
strand $c_{\tilde k}$ or strand $b$, which is impossible, since
$w(\tilde R)$ is an $r$-semi-shuffle. Hence (\ref{eq:check3})
holds.

If $c_{\tilde k-1}\leq r$, then strands $c=c_{\tilde k}$ and $d=b$
pass next to each other in row $\ell+1$ at a place $(\ell,j)$, so
insertion into $(\ell,j)$ is allowed. It implies $b_{i^*}\geq b$.

It remains to consider the case when $c_{\tilde k-1}=b-1$. Let
$\bar b_k$ denote the entries of $B(T(\ell+1))$. We will prove
that
\begin{equation}
\label{eq:check11} i_*\neq 1\text{ and }c_{\tilde k-1}=\bar
b_{i-1}=b-1.\end{equation} If (\ref{eq:check11}) holds, then,
since step $(\ell,i_*)$ is an insertion, $i_*$ is not the first
box in row $\ell$ and box $i_*-1$ is in the same row as box $i_*$.
Hence by induction assumption
$$
b-1=\bar b_{i-1}\leq b_{i_*-1}<b_{i_*}.
$$
Therefore, since $b> \bar b_{i-1}=b-1$, we conclude $b=b_{i}\leq
b_{i_*}$.

It remain to show that if $c_{\tilde k-1}=b-1$ then
(\ref{eq:check11}) holds. Since $w(\tilde R)$ is an
$r$-semi-shuffle and strands $b$ and $b-1$ pass next to each other
in row $\ell$ of $\tilde R$,
$$w(\tilde R_{\geq\ell+1})(b)-1=w(\tilde R_{\geq\ell+1})(b-1).$$

$R(\ell,i-1)_{\geq\ell+1}$ is constructed out of $\tilde
R_{\geq\ell+1}$ by adding some crossings. It is not difficult to
see that if (\ref{eq:check11}) fails, none of these crossings
involve strands $b$ or $b-1$.~So
$$w(R(\ell,i-1)_{\geq\ell+1})(b)-1=w(R(\ell,i-1)_{\geq\ell+1})(b-1).$$
But it is impossible by (\ref{no i}), since $a<b-1<b$ and
$$
w(R(\ell,i-1)_{\geq\ell+1})(a)<w(R(\ell,i-1)_{\geq\ell+1})(b-1)<w(R(\ell,i-1)_{\geq\ell+1})(b)
$$
while $l(w(R(\ell,i-1)_{\geq\ell+1})t_{ab})=
l(w(R(\ell,i-1)_{\geq\ell+1}))+1$.

\subsubsection{Case (\ref{case1}).}
\label{proof(4),case2} As before, let $(ab)$ and $(ef)$ be the
entries of boxes $i$ and $i-1$ of~$T'$. Assume $(\ell,i-1)$ is a
rectification (the argument below can be easily modified to
provide a proof in the case step $(\ell,i-1)$ is an insertion).
Assume $(\bar e, \bar f)$ is the entry of box $i_+-1$ of
$T(\ell+1)$.

Start with the case when box $i$ is not in the first column of
$sh(T')$. Then let us show that if $R'$ is an rc-graph, then
$E(w(R_{\geq\ell}),T')$ is row and column strict. It is obvious if
$(ef)\neq (\bar e \bar f)$, since in this case the entry of box
$i-1$ of $E(w(R_{\geq\ell}),T')$ is smaller then the value of box
$i_+-1$ of $E(w(R_{\geq\ell}),T(\ell+1))$.

If $(ef)=(\bar e\bar f)$ and $E(w(R_{\geq\ell}),T')$ is not row
and column strict, then, strands $b$ and $f$ intersect in row
$\ell$ in $R'$, such that $f \boxplus b=\ell$. But then $f\boxplus
a$ in $R(\ell,i-1)$, which is impossible since $a<f$.

In the case when $i$ is in the first column and $R'$ is an
rc-graph, we will show one of the following holds:
\begin{align}
\label{eq:check4-1}&\text{$w(R(\ell,i-1))(f)>w(R(\ell,i-1))( a)$,}\\
\label{eq:check4-2}&\text{$a=e$ and $a\boxplus f=\ell$ in
$R(\ell,i-1)$,}
\end{align}
If (\ref{eq:check4-1}) holds, then $E(w(R_{\geq\ell}),T')$ is row
and column strict. If (\ref{eq:check4-2}) holds, then $b\boxplus
f=\ell$ in $R'$ as shown in the second picture of Figure
\ref{remove-crossing}. So, it remain to prove (\ref{eq:check4-1})
or (\ref{eq:check4-2}) hold in the case when $i$ is in the first
column and $R'$ is an rc-graph.

Since $E(w(R_{\geq\ell+1}), T(\ell+1))$ is a row and column strict
tableau, we know
\begin{equation}
\label{eq:inequality} w(S_{i_+-1}(\mathcal P(\ell+1)))(\bar
f)>w(S_{i_+}(\mathcal P(\ell+1)))(b) =w(S_{i_+-1}(\mathcal
P(\ell+1)))( a).
\end{equation}
Moreover, removing crossings from row $\ell$ of $R(\ell,i-1)$
produces $S_{i_+-1}(\mathcal P(\ell+1))$.

In the case $(ef)=(\bar e\bar f)$, the inequality
(\ref{eq:inequality}) implies (\ref{eq:check4-1}), unless
$a\boxplus f=\ell$ in $R(\ell,i-1)$. But it is not difficult to
see that if $(ef)=(\bar e\bar f)$, strands $a$ and $f$ cannot
intersect.

Otherwise, if $(ef)\neq(\bar e\bar f)$, we will show during the
proof of (5) that a crossing of strands $\bar e,\bar f$ has been
removed during step $(\ell,i-1)$ from place $(\ell,\tilde j)$ and
then another crossing has been inserted to the left of
$(\ell,\tilde j)$. Assume $\bar R'$ is the intermediate rc-graph
in step $(\ell,i-1)$ constructed by removing a crossing from
$R(\ell,i-2)$. Consider sequence $c_k$ defined by (\ref{sequence
ck}) for $\bar R'$. Let $c_0=a$, if $a$ does not intersect row
$\ell$ vertically, otherwise set $c_0=a'$ with $a'\boxplus
a=\ell$. Then $\bar f$ is an element of the sequence $c_k$. Let
$\bar f=c_{\tilde k}$. Clearly, there exist a place $(\ell,j)$,
where insertion is allowed with strands $c_{\bar k}$ and $c_{\bar
k+1}$ passing through $(\ell,j)$ with $0\leq\bar k<\tilde k$.
Choose such place with the largest possible $j$, let it be
$(\ell,j_1)$, then we define $R(\ell,i-1)$ and $T(\ell,i-1)$ is
such a way that $e=c_{\bar k}$ and $f=c_{\bar k+1}$ . If $\bar
k=0$ and $c_0=a$, (\ref{eq:check4-2}) holds, otherwise
(\ref{eq:check4-1}) must be satisfied.

\subsection{Proof of (5)}
Assume that a crossing at place $ (\ell,j_0)$ in $R'$ has been
removed to produce $R''$. It is not difficult to see that $R''$ is
an rc-graph, $(w(R_{\geq\ell}),R'',T'')$ is an $r$-Bruhat package
and $E(w(R_{\geq\ell}), T')$ is row and column strict. We need to
show that there exist a place where insertion is allowed to the
left of $(\ell,j_0)$ and after $R(\ell,i)$ and $T(\ell,i)$ are
defined, (\ref{main-cond}) and (\ref{cond1}) are satisfied.

\subsubsection{Case (\ref{case2})}
As in the proof of (3), we can use sequence $c_k$ for $R''$ to
show that a place, where insertion is allowed, to the left of
place $(\ell,j_0)$ exists. Moreover, the rightmost place
$(\ell,j_1)$ where insertion is allowed looks like the first
picture in Figure \ref{allowed-insertions-rect}.

After inserting crossing $(\ell,j_1)$ into $R''$ to defining
$R(\ell,i)$ and $(cd)$ into box $i$ of $T''$ to define
$T(\ell,i-1)$, it is easy to see $R(\ell,i)$ is an rc-graph and
(\ref{cond1}) is satisfied. By Lemma \ref{lemma:standard} to show
(\ref{main-cond}) holds, it is enough to show $B(T(\ell,i))$ is
row strict. Notice that $B(T(\ell,i))$ differs from
$B(T(\ell,i-1))$ only in box $i$. So we just have to check that
the entry box $i$ is still greater than the entry of the box to
the left of box $i$ and not greater than the entry of the box
below box $i$. This can be done by an argument, which is almost
identical to the argument used in Section \ref{proof(4),case1}.

\subsubsection{Case (\ref{case1})}
Recall that $(ab)$ is the entry of box $i$ of $T'$, $(eg)$ is the
entry of box $i-1$ of $T''$. Consider sequence $c_k$ for $R''$,
defined by (\ref{sequence ck}), which starts with $c_0=\ell$ and
ends with $c_{\tilde k}=b$. Then there exists $\bar k$ between $0$
and $\tilde k-1$, such that strand $c_{\bar k}$ and $c_{\bar k+1}$
pass next to each other in a place where insertion is allowed. Let
the rightmost place to the left of $(\ell,j_0)$, where insertion
is allowed, be $(\ell,j_1)$.

Consider the case when box $i$ is in the first column of $T''$.
Then using sequence $c_k$, it is easy to see that place
$(\ell,j_1)$ looks like the first pictures from Figure
\ref{allowed-insertions-rect}. (We used this in Section
\ref{proof(4),case2}.) Therefore, as for the insertion step,
$R(\ell,i)$ is an rc-graph and (\ref{cond1}) holds. Moreover
(\ref{main-cond}) holds, since strand $g$ passes row $\ell$ to the
right of place $(\ell,j_1)$.

Otherwise, if box $i$ is in the first row, but not the first
element of this row, then strand $g$ is either an element of the
sequence $c_k$ or intersects one of the strand $c_k$ in row
$\ell$. So place $(\ell,j_1)$ could look like the first picture of
Figure \ref{allowed-insertions-rect} and strand $g$ passes to the
left of this place. Or, it could look like the second picture of
Figure \ref{allowed-insertions-rect}.

If it is the first picture, then, as before, $R''$ is an rc-graph,
(\ref{cond1}) holds, while (\ref{main-cond}) holds, since strand
$g$ passes row $\ell$ to the left of place $(\ell,j_1)$.

If it is the second picture, it is easy to see that $R''$ is an
rc-graph and that (\ref{main-cond}) holds, while to prove
(\ref{cond1}), we must show
\begin{equation}
\label{check10} l(w(R(\ell,i)))=l(w(R(\ell,i))t_{eg}t_{ ed})+2
=l(w(R(\ell,i))t_{e g })+1.
\end{equation}
To prove the first equality of (\ref{check10}), notice
$t_{eg}t_{ed}=t_{gd}t_{eg}$, hence
\begin{align*}
l(w(R(\ell,i))t_{eg}t_{ ed})&=l(w(R(\ell,i))t_{gd}t_{eg})
=l(w(R'')t_{eg})\\ &=l(w(R''))-1= l(w(R(\ell,i)))-2.
\end{align*}
For the second equality, notice that $e<g$ and
$w(R(\ell,i))(e)>w(R(\ell,i))(g)$, thus
$$
l(w(R(\ell,i))t_{e g })<l(w(R(\ell,i))).
$$
At the same time, $e<d$ and
$w(R(\ell,i))t_{eg}(e)>w(R(\ell,i))t_{eg}(d)$, hence
$$
l(w(R(\ell,i)t_{e f }t_{ed})<l(w(R(\ell,i))t_{eg}).
$$
This proves the second part of (\ref{check10}).

\section{Inverse Insertion Algorithm}
\label{inverse algorithm} Given an rc-graph $U$ and a tableau of
transposition $T$ of $v(\lambda,r)$ inverse insertion algorithm
defines rc-graphs $U\rightarrow T$ and $Y(R,T)$, given that $T$ is
an $r$-Brihat chain of $w=w(U)w(T)^{-1}$, $E(w,T)$ is row and
column strict, and $w,v(\lambda,r)$ satisfy (\ref{case2}) or
(\ref{case1}). This section describes the inverse algorithm.

\subsection{Sequence of inverse steps} Inverse insertion algorithm
performs the same steps as insertion algorithm but in the opposite
order. Each step will be either an inverse row-to-to step, an
inverse insertion or an inverse rectification.

Each step $(\ell,i)$ with $1\leq i\leq m_\ell$ constructs rc-graph
$R(\ell,i-1)$ with no crossings above row $\ell$ and tableau of
transposition $T(\ell,i-1)$ filled up to $i-1$. Each
step~$(\ell,0)$ defines integer $m_{\ell+1}$, an rc-graph
$R(\ell+1,m_{\ell+1})$ with no crossings above row $\ell+1$ and a
tableau of transpositions $T(\ell+1,m_{\ell+1})$. Conditions
(\ref{main-cond}), (\ref{cond1}) always hold.

\subsection{Start of the algorithm} Set $m_1=|T|$, $R(1,m_1)=U$ and
$T(1,m_1)=T$.

\subsection{Inverse insertion}
\label{reverse insertion} 
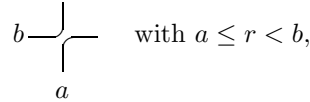
\begin{figure}[ht]
\begin{picture}(100, 25)

\put (0,15){\line(1,0){10}} \put (13,2){\line(0,1){10}} \put
(16,15){\line(1,0){10}} \put (13,18){\line(0,1){10}}

\put(9,19){\oval(8,8)[br]} \put(17,11){\oval(8,8)[tl]}

\put(-6,13){$b$} \put(10,-8){$a$}

\put(40,13){with $a\leq r <b$,}

\end{picture}
\caption{Place where insertion is allowed} \label{inverse-allowed}
\end{figure}

Consider step $(\ell,i)$ with $i>0$. We
need to construct $R(\ell,i-1)$ and $T(\ell,i-1)$. Let $(cd)$ be
the entry of box $i$ of $T(\ell,i)$. By Lemma~\ref{add-crossing},
strand $c$ and $d$ intersect in $R(\ell,i)$ at some place
$(\ell_0,j_0)$. If $\ell=\ell_0$ define $T''$ by removing the
entry of box $i$ from $T(\ell,i)$. Define $R''$ by removing the
crossing of strands $c$ and $d$ from $R(\ell,i)$ from place
$(\ell,j_0)$. We say that insertion into place $(\ell,j)$ is
allowed if strands $a,b$ pass this place as shown in Figure
\ref{inverse-allowed}.
If there are no places $(\ell,j)$, where insertion is allowed,
with $j>j_0$, this step is an inverse insertion, which sets
$R(\ell,i-1)=R''$ and $T(\ell,i-1)=T''$.

\subsection{Inverse rectification}
All steps $(\ell,i)$ with $i>0$, which are not inverse insertions
are inverse rectifications.

Adopt the notation from previous section. If $\ell_0\neq\ell$,
define $T(\ell,i-1)$ by emptying box $i$ of $T(\ell,i)$ and define
$R(\ell,i-1)$ by removing the crossing of strands $c$ and $d$ and
move on to the next step, except for one case. Namely, if $(ef)$
is the entry of $i-1$ of $T(\ell,i)$, $c=e$ and $f\boxplus
d=\ell$. In this case define $R''$ by removing the crossing of $b$
and $f$ and define $T''$ by emptying box $i$ of $T(\ell,i)$ and
placing $(ed)=(cd)$ in box $i-1$. If $\ell_0=\ell$, define $R'$
and $T'$ as it was done in the previous section.

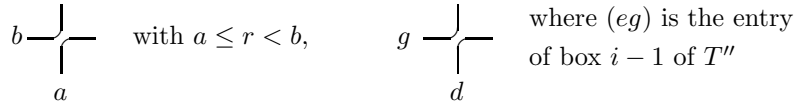
\begin{figure}[ht]
\begin{picture}(320, 25)

\put (0,15){\line(1,0){10}} \put (13,2){\line(0,1){10}} \put
(16,15){\line(1,0){10}} \put (13,18){\line(0,1){10}}

\put(9,19){\oval(8,8)[br]} \put(17,11){\oval(8,8)[tl]}

\put(-6,13){$b$} \put(10,-8){$a$}

\put(40,13){with $a\leq r <b$,}


\put (150,15){\line(1,0){10}} \put (163,2){\line(0,1){10}} \put
(166,15){\line(1,0){10}} \put (163,18){\line(0,1){10}}

\put(159,19){\oval(8,8)[br]} \put(167,11){\oval(8,8)[tl]}

\put(140,13){$g$} \put(160,-8){$d$}

\put(190,20){where $(eg)$ is the entry} \put(190,5){of box $i-1$
of $T''$ }

\end{picture}
\caption{Places, where insertions are allowed during
rectification.} \label{inverse:allowed-insertions-rect}
\end{figure}

Once $R''$ and $T''$ are constructed, we say that insertion into
places in row $\ell$ shown in
Figure~\ref{inverse:allowed-insertions-rect} are allowed.
Find the leftmost place where insertion is allowed to the right of
place $(\ell,j_0)$. Insert a crossing into this place to define
$R'$. If this place looks like the first picture of Figure
\ref{inverse:allowed-insertions-rect}, add $(ab)$ to box $i$ of
$T''$ to construct $T'$, otherwise insert $(ed)$ and $(eg)$ into
boxes $i-1$ and $i$ of $T''$ to produce $T'$.

Once $R'$ and $T'$ are constructed, let $(ab)$ be the entry of box
$i$ of $T'$. Then it can be shown that strands $a$ and $b$
intersect below row $\ell$.  Remove this crossing to produce
$R(\ell,i-1)$ and construct $T(\ell,i-1)$ by emptying box $i$ of
$T'$.

\subsection{Inverse row-to-row steps} Steps $(\ell,0)$ are
inverse row-to-row steps. They define $m_\ell$ to be the number of
inverse rectifications $(\ell,i_1),\dots,(\ell,i_{m_\ell+1})$ for
row $\ell$. Also each step $(\ell,0)$ sets
$R(\ell+1,m_{\ell+1})=R(\ell,m_\ell)_{\geq\ell+1}$.

The shape  of $T(\ell+1,m_{\ell+1})$ is the subdiagram of
$sh(T(\ell,m_\ell))$ consisting of boxes $(i_1,\dots,
i_{m_{\ell+1}})$. By construction, this will be a Young diagram.
The entry $(a_{k}b_k)$ of box $k$ of $T(\ell+1,m_{\ell+1})$ is
determined by
$$
w(R(\ell,i_{k-1})_{\geq\ell+1})=w(R(\ell,i_{k})_{\geq\ell+1})t_{a_kb_k}.
$$

\subsection{End of the inverse algorithm} Set $U\rightarrow
T=R(r,0)$. We will define $Y(U,T)$ by presenting its word. Set
$word_{r+1}$ to be empty. Define $word_\ell$ of length $m_\ell$ by
adding letters $\ell$ to $word_{\ell+1}$ as follows. If
$(\ell,i_1),\dots,(\ell,i_{m_\ell+1})$ are the rectification steps
for  row $\ell$. Then set letter $i_k$ of $word_\ell$ to be the
same as letter $k$ of $word_{\ell+1}$, set all the other letters
of $word_\ell$ to be equal to $\ell$. Finally set
$word(Y(U,T))=word_1$.

\section{Examples}
\label{examples}

\subsection{Example of rc-graphs $S_j(\mathcal P)$}
\label{example1} 
\begin{figure}[ht]
\begin{picture}(200,88)


\put(0,64){1} \put(0,48){2} \put(0,32){3} \put(0,16){4}
\put(0,0){5}

\put(18,83){1} \put(34,83){2} \put(50,83){3} \put(66,83){4}
\put(82,83){5}

\put(81,71){\oval(8,8)[br]} \put(65,55){\oval(8,8)[br]}
\put(49,39){\oval(8,8)[br]} \put(33,23){\oval(8,8)[br]}
\put(17,7){\oval(8,8)[br]}

\put(65,71){\oval(8,8)[br]} \put(73,63){\oval(8,8)[tl]}
\put(33,39){\oval(8,8)[br]} \put(41,31){\oval(8,8)[tl]}
\put(17,23){\oval(8,8)[br]} \put(25,15){\oval(8,8)[tl]}
\put(17,39){\oval(8,8)[br]} \put(25,31){\oval(8,8)[tl]}


\put (18,51){\line(1,0){6}} \put (21,48){\line(0,1){6}}

\put (18,67){\line(1,0){6}} \put (21,64){\line(0,1){6}}

\put (34,51){\line(1,0){6}} \put (37,48){\line(0,1){6}}

\put (34,67){\line(1,0){6}} \put (37,64){\line(0,1){6}}

\put (50,51){\line(1,0){6}} \put (53,48){\line(0,1){6}}

\put (50,67){\line(1,0){6}} \put (53,64){\line(0,1){6}}


\put (08,67){\line(1,0){10}} \put (08,51){\line(1,0){10}} \put
(08,35){\line(1,0){10}} \put (08,19){\line(1,0){10}} \put
(08,3){\line(1,0){10}} \put (24,67){\line(1,0){10}} \put
(24,51){\line(1,0){10}} \put (24,35){\line(1,0){10}} \put
(24,19){\line(1,0){10}} \put (40,67){\line(1,0){10}} \put
(40,51){\line(1,0){10}} \put (40,35){\line(1,0){10}} \put
(56,67){\line(1,0){10}} \put (56,51){\line(1,0){10}} \put
(72,67){\line(1,0){10}}


\put (21,70){\line(0,1){10}} \put (21,54){\line(0,1){10}} \put
(21,38){\line(0,1){10}} \put (21,22){\line(0,1){10}} \put
(21,6){\line(0,1){10}} \put (37,70){\line(0,1){10}} \put
(37,54){\line(0,1){10}} \put (37,38){\line(0,1){10}} \put
(37,22){\line(0,1){10}} \put (53,70){\line(0,1){10}} \put
(53,54){\line(0,1){10}} \put (53,38){\line(0,1){10}} \put
(69,70){\line(0,1){10}} \put (69,54){\line(0,1){10}} \put
(85,70){\line(0,1){10}}

\put (130,55){\line(1,0){15}} \put (130,40){\line(1,0){15}} \put
(145,55){\line(1,0){15}} \put (145,40){\line(1,0){15}} \put
(130,25){\line(1,0){15}} \put (145,25){\line(1,0){15}} \put
(145,25){\line(1,0){15}}

\put (130,40){\line(0,1){15}} \put (145,40){\line(0,1){15}} \put
(160,40){\line(0,1){15}} \put (130,25){\line(0,1){15}} \put
(145,25){\line(0,1){15}} \put (160,25){\line(0,1){15}}

\put(133,29){14} \put(148,29){23} \put(133,44){25}
\put(148,44){15}

\end{picture}
\caption{Rc-graph $R$ and tableau of transpositions $T$.}
\label{example1-figure1}
\end{figure}
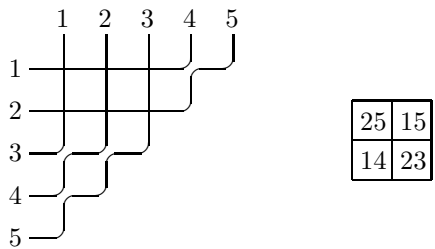

Assume $R$ and $T$ are given in
Figure~\ref{example1-figure1}. Define
$w=w(R)w(T)^{-1}=(2,1,4,3,5,6,\dots)$. Then $\mathcal P=(w,R,T)$
is an $r$-Bruhat package.

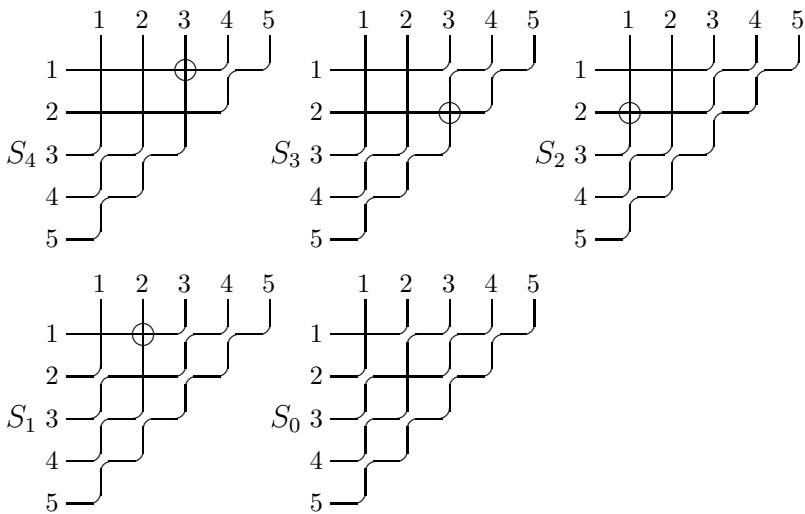
\begin{figure}[ht]
\begin{picture}(300,188)

\put(-15,132){\Large $S_4$} \put(85,132){\Large $S_3$}
\put(185,132){\Large $S_2$} \put(-15,32){\Large $S_1$}
\put(85,32){\Large $S_0$}


\put(0,164){1} \put(0,148){2} \put(0,132){3} \put(0,116){4}
\put(0,100){5}

\put(18,183){1} \put(34,183){2} \put(50,183){3} \put(66,183){4}
\put(82,183){5}

\put(81,171){\oval(8,8)[br]} \put(65,155){\oval(8,8)[br]}
\put(49,139){\oval(8,8)[br]} \put(33,123){\oval(8,8)[br]}
\put(17,107){\oval(8,8)[br]}

\put(65,171){\oval(8,8)[br]} \put(73,163){\oval(8,8)[tl]}
\put(33,139){\oval(8,8)[br]} \put(41,131){\oval(8,8)[tl]}
\put(17,123){\oval(8,8)[br]} \put(25,115){\oval(8,8)[tl]}
\put(17,139){\oval(8,8)[br]} \put(25,131){\oval(8,8)[tl]}


\put (18,151){\line(1,0){6}} \put (21,148){\line(0,1){6}}

\put (18,167){\line(1,0){6}} \put (21,164){\line(0,1){6}}

\put (34,151){\line(1,0){6}} \put (37,148){\line(0,1){6}}

\put (34,167){\line(1,0){6}} \put (37,164){\line(0,1){6}}

\put (50,151){\line(1,0){6}} \put (53,148){\line(0,1){6}}

\put (50,167){\line(1,0){6}} \put (53,164){\line(0,1){6}}


\put (08,167){\line(1,0){10}} \put (08,151){\line(1,0){10}} \put
(08,135){\line(1,0){10}} \put (08,119){\line(1,0){10}} \put
(08,103){\line(1,0){10}} \put (24,167){\line(1,0){10}} \put
(24,151){\line(1,0){10}} \put (24,135){\line(1,0){10}} \put
(24,119){\line(1,0){10}} \put (40,167){\line(1,0){10}} \put
(40,151){\line(1,0){10}} \put (40,135){\line(1,0){10}} \put
(56,167){\line(1,0){10}} \put (56,151){\line(1,0){10}} \put
(72,167){\line(1,0){10}}


\put (21,170){\line(0,1){10}} \put (21,154){\line(0,1){10}} \put
(21,138){\line(0,1){10}} \put (21,122){\line(0,1){10}} \put
(21,106){\line(0,1){10}} \put (37,170){\line(0,1){10}} \put
(37,154){\line(0,1){10}} \put (37,138){\line(0,1){10}} \put
(37,122){\line(0,1){10}} \put (53,170){\line(0,1){10}} \put
(53,154){\line(0,1){10}} \put (53,138){\line(0,1){10}} \put
(69,170){\line(0,1){10}} \put (69,154){\line(0,1){10}} \put
(85,170){\line(0,1){10}}

\put(53,167){\oval(8,8)}


\put(100,164){1} \put(100,148){2} \put(100,132){3}
\put(100,116){4} \put(100,100){5}

\put(118,183){1} \put(134,183){2} \put(150,183){3}
\put(166,183){4} \put(182,183){5}

\put(181,171){\oval(8,8)[br]} \put(165,155){\oval(8,8)[br]}
\put(149,139){\oval(8,8)[br]} \put(133,123){\oval(8,8)[br]}
\put(117,107){\oval(8,8)[br]}

\put(165,171){\oval(8,8)[br]} \put(173,163){\oval(8,8)[tl]}
\put(133,139){\oval(8,8)[br]} \put(141,131){\oval(8,8)[tl]}
\put(117,123){\oval(8,8)[br]} \put(125,115){\oval(8,8)[tl]}
\put(117,139){\oval(8,8)[br]} \put(125,131){\oval(8,8)[tl]}
\put(149,171){\oval(8,8)[br]} \put(157,163){\oval(8,8)[tl]}


\put (118,151){\line(1,0){6}} \put (121,148){\line(0,1){6}}

\put (118,167){\line(1,0){6}} \put (121,164){\line(0,1){6}}

\put (134,151){\line(1,0){6}} \put (137,148){\line(0,1){6}}

\put (134,167){\line(1,0){6}} \put (137,164){\line(0,1){6}}

\put (150,151){\line(1,0){6}} \put (153,148){\line(0,1){6}}


\put (108,167){\line(1,0){10}} \put (108,151){\line(1,0){10}} \put
(108,135){\line(1,0){10}} \put (108,119){\line(1,0){10}} \put
(108,103){\line(1,0){10}} \put (124,167){\line(1,0){10}} \put
(124,151){\line(1,0){10}} \put (124,135){\line(1,0){10}} \put
(124,119){\line(1,0){10}} \put (140,167){\line(1,0){10}} \put
(140,151){\line(1,0){10}} \put (140,135){\line(1,0){10}} \put
(156,167){\line(1,0){10}} \put (156,151){\line(1,0){10}} \put
(172,167){\line(1,0){10}}


\put (121,170){\line(0,1){10}} \put (121,154){\line(0,1){10}} \put
(121,138){\line(0,1){10}} \put (121,122){\line(0,1){10}} \put
(121,106){\line(0,1){10}} \put (137,170){\line(0,1){10}} \put
(137,154){\line(0,1){10}} \put (137,138){\line(0,1){10}} \put
(137,122){\line(0,1){10}} \put (153,170){\line(0,1){10}} \put
(153,154){\line(0,1){10}} \put (153,138){\line(0,1){10}} \put
(169,170){\line(0,1){10}} \put (169,154){\line(0,1){10}} \put
(185,170){\line(0,1){10}}

\put(153,151){\oval(8,8)}

\put(200,164){1} \put(200,148){2} \put(200,132){3}
\put(200,116){4} \put(200,100){5}

\put(218,183){1} \put(234,183){2} \put(250,183){3}
\put(266,183){4} \put(282,183){5}

\put(281,171){\oval(8,8)[br]} \put(265,155){\oval(8,8)[br]}
\put(249,139){\oval(8,8)[br]} \put(233,123){\oval(8,8)[br]}
\put(217,107){\oval(8,8)[br]}

\put(265,171){\oval(8,8)[br]} \put(273,163){\oval(8,8)[tl]}
\put(233,139){\oval(8,8)[br]} \put(241,131){\oval(8,8)[tl]}
\put(217,123){\oval(8,8)[br]} \put(225,115){\oval(8,8)[tl]}
\put(217,139){\oval(8,8)[br]} \put(225,131){\oval(8,8)[tl]}
\put(249,171){\oval(8,8)[br]} \put(257,163){\oval(8,8)[tl]}
\put(249,155){\oval(8,8)[br]} \put(257,147){\oval(8,8)[tl]}


\put (218,167){\line(1,0){6}} \put (221,164){\line(0,1){6}}

\put (234,151){\line(1,0){6}} \put (237,148){\line(0,1){6}}

\put (234,167){\line(1,0){6}} \put (237,164){\line(0,1){6}}

\put (218,151){\line(1,0){6}} \put (221,148){\line(0,1){6}}


\put (208,167){\line(1,0){10}} \put (208,151){\line(1,0){10}} \put
(208,135){\line(1,0){10}} \put (208,119){\line(1,0){10}} \put
(208,103){\line(1,0){10}} \put (224,167){\line(1,0){10}} \put
(224,151){\line(1,0){10}} \put (224,135){\line(1,0){10}} \put
(224,119){\line(1,0){10}} \put (240,167){\line(1,0){10}} \put
(240,151){\line(1,0){10}} \put (240,135){\line(1,0){10}} \put
(256,167){\line(1,0){10}} \put (256,151){\line(1,0){10}} \put
(272,167){\line(1,0){10}}


\put (221,170){\line(0,1){10}} \put (221,154){\line(0,1){10}} \put
(221,138){\line(0,1){10}} \put (221,122){\line(0,1){10}} \put
(221,106){\line(0,1){10}} \put (237,170){\line(0,1){10}} \put
(237,154){\line(0,1){10}} \put (237,138){\line(0,1){10}} \put
(237,122){\line(0,1){10}} \put (253,170){\line(0,1){10}} \put
(253,154){\line(0,1){10}} \put (253,138){\line(0,1){10}} \put
(269,170){\line(0,1){10}} \put (269,154){\line(0,1){10}} \put
(285,170){\line(0,1){10}}

\put(221,151){\oval(8,8)}


\put(0,64){1} \put(0,48){2} \put(0,32){3} \put(0,16){4}
\put(0,0){5}

\put(18,83){1} \put(34,83){2} \put(50,83){3} \put(66,83){4}
\put(82,83){5}

\put(81,71){\oval(8,8)[br]} \put(65,55){\oval(8,8)[br]}
\put(49,39){\oval(8,8)[br]} \put(33,23){\oval(8,8)[br]}
\put(17,7){\oval(8,8)[br]}

\put(65,71){\oval(8,8)[br]} \put(73,63){\oval(8,8)[tl]}
\put(33,39){\oval(8,8)[br]} \put(41,31){\oval(8,8)[tl]}
\put(17,23){\oval(8,8)[br]} \put(25,15){\oval(8,8)[tl]}
\put(17,39){\oval(8,8)[br]} \put(25,31){\oval(8,8)[tl]}
\put(49,71){\oval(8,8)[br]} \put(57,63){\oval(8,8)[tl]}
\put(17,55){\oval(8,8)[br]} \put(25,47){\oval(8,8)[tl]}
\put(49,55){\oval(8,8)[br]} \put(57,47){\oval(8,8)[tl]}


\put (18,67){\line(1,0){6}} \put (21,64){\line(0,1){6}}

\put (34,51){\line(1,0){6}} \put (37,48){\line(0,1){6}}

\put (34,67){\line(1,0){6}} \put (37,64){\line(0,1){6}}


\put (08,67){\line(1,0){10}} \put (08,51){\line(1,0){10}} \put
(08,35){\line(1,0){10}} \put (08,19){\line(1,0){10}} \put
(08,3){\line(1,0){10}} \put (24,67){\line(1,0){10}} \put
(24,51){\line(1,0){10}} \put (24,35){\line(1,0){10}} \put
(24,19){\line(1,0){10}} \put (40,67){\line(1,0){10}} \put
(40,51){\line(1,0){10}} \put (40,35){\line(1,0){10}} \put
(56,67){\line(1,0){10}} \put (56,51){\line(1,0){10}} \put
(72,67){\line(1,0){10}}


\put (21,70){\line(0,1){10}} \put (21,54){\line(0,1){10}} \put
(21,38){\line(0,1){10}} \put (21,22){\line(0,1){10}} \put
(21,6){\line(0,1){10}} \put (37,70){\line(0,1){10}} \put
(37,54){\line(0,1){10}} \put (37,38){\line(0,1){10}} \put
(37,22){\line(0,1){10}} \put (53,70){\line(0,1){10}} \put
(53,54){\line(0,1){10}} \put (53,38){\line(0,1){10}} \put
(69,70){\line(0,1){10}} \put (69,54){\line(0,1){10}} \put
(85,70){\line(0,1){10}}
\put(37,67){\oval(8,8)}

\put(100,64){1} \put(100,48){2} \put(100,32){3} \put(100,16){4}
\put(100,0){5}

\put(118,83){1} \put(134,83){2} \put(150,83){3} \put(166,83){4}
\put(182,83){5}

\put(181,71){\oval(8,8)[br]} \put(165,55){\oval(8,8)[br]}
\put(149,39){\oval(8,8)[br]} \put(133,23){\oval(8,8)[br]}
\put(117,7){\oval(8,8)[br]}

\put(165,71){\oval(8,8)[br]} \put(173,63){\oval(8,8)[tl]}
\put(133,39){\oval(8,8)[br]} \put(141,31){\oval(8,8)[tl]}
\put(117,23){\oval(8,8)[br]} \put(125,15){\oval(8,8)[tl]}
\put(117,39){\oval(8,8)[br]} \put(125,31){\oval(8,8)[tl]}
\put(149,71){\oval(8,8)[br]} \put(157,63){\oval(8,8)[tl]}
\put(133,71){\oval(8,8)[br]} \put(141,63){\oval(8,8)[tl]}
\put(149,55){\oval(8,8)[br]} \put(157,47){\oval(8,8)[tl]}
\put(117,55){\oval(8,8)[br]} \put(125,47){\oval(8,8)[tl]}

\put (118,67){\line(1,0){6}} \put (121,64){\line(0,1){6}}

\put (134,51){\line(1,0){6}} \put (137,48){\line(0,1){6}}


\put (108,67){\line(1,0){10}} \put (108,51){\line(1,0){10}} \put
(108,35){\line(1,0){10}} \put (108,19){\line(1,0){10}} \put
(108,3){\line(1,0){10}} \put (124,67){\line(1,0){10}} \put
(124,51){\line(1,0){10}} \put (124,35){\line(1,0){10}} \put
(124,19){\line(1,0){10}} \put (140,67){\line(1,0){10}} \put
(140,51){\line(1,0){10}} \put (140,35){\line(1,0){10}} \put
(156,67){\line(1,0){10}} \put (156,51){\line(1,0){10}} \put
(172,67){\line(1,0){10}}


\put (121,70){\line(0,1){10}} \put (121,54){\line(0,1){10}} \put
(121,38){\line(0,1){10}} \put (121,22){\line(0,1){10}} \put
(121,6){\line(0,1){10}} \put (137,70){\line(0,1){10}} \put
(137,54){\line(0,1){10}} \put (137,38){\line(0,1){10}} \put
(137,22){\line(0,1){10}} \put (153,70){\line(0,1){10}} \put
(153,54){\line(0,1){10}} \put (153,38){\line(0,1){10}} \put
(169,70){\line(0,1){10}} \put (169,54){\line(0,1){10}} \put
(185,70){\line(0,1){10}}

\end{picture}
\caption{Rc-graph $S_4(\mathcal P)$ through $S_0(\mathcal P)$.}
\label{example1-figure2}
\end{figure}

Then sequence $S_{j}(\mathcal P)$ is given in Figure
\ref{example1-figure2}. In each graph $S_j(\mathcal P)$ the
circled  crossing needs to be removed to construct
$S_{j-1}(\mathcal P)$. Since $word(T)=(14)(23)(25)(15)$, $S_3$ is
constructed out of $S_4$ by removing the crossing of strands $1$
and $5$, $S_2$ out of $S_3$ by removing the crossing of strands
$2$ and $5$, and so on.

\newpage
\subsection{Example of insertion algorithm in case (\ref{case2})}
\label{example2} From now on we draw only crossings of rc-graphs
without drawing strands, as it was done in Figure~\ref{definition
of graphs}. It makes it easier to see how rc-graphs change during
the algorithm. At the same time, as usual, we assume each rc-graph
extends infinitely to the right and to the bottom and the part of
each rc-graph, which is not shown, has no crossings.

Assume the rc-graphs $R$ and $Y$ are given in Figure
\ref{ex2:R,T}, so that $r=3$, $w(R)=(1,4,3,2,5,6,\dots)$ is a
$3$-semi-shuffle and $w(Y)=(1,4,5,2,3,6,\dots)=v((2,2),3)$ is a
$3$-shuffle. We will to illustrate all the steps of the algorithm
for $R\leftarrow Y$.

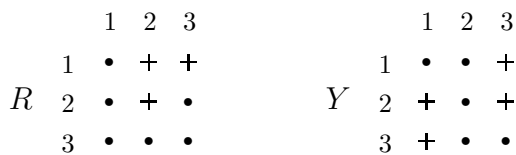
\begin{figure}[ht]
\begin{picture}(160,56)

\put(0,29){1} \put(0,14){2} \put(0,-1){3}

\put(16,45){1} \put(31,45){2} \put(46,45){3}


\put (30,18){\line(1,0){6}} \put (33,15){\line(0,1){6}}

\put (30,33){\line(1,0){6}} \put (33,30){\line(0,1){6}}

\put (45,33){\line(1,0){6}} \put (48,30){\line(0,1){6}}


\put(18,33){\circle*{3}} \put(18,3){\circle*{3}}
\put(48,18){\circle*{3}} \put(18,18){\circle*{3}}
\put(33,3){\circle*{3}}\put(48,3){\circle*{3}}

\put(-20,15){\Large $R$}



\put(120,29){1} \put(120,14){2} \put(120,-1){3}

\put(136,45){1} \put(151,45){2} \put(166,45){3}


\put (135,3){\line(1,0){6}} \put (138,0){\line(0,1){6}}

\put (135,18){\line(1,0){6}} \put (138,15){\line(0,1){6}}

\put (165,33){\line(1,0){6}} \put (168,30){\line(0,1){6}}

\put (165,18){\line(1,0){6}} \put (168,15){\line(0,1){6}}


\put(138,33){\circle*{3}} \put(153,18){\circle*{3}}
\put(153,33){\circle*{3}} \put(153,3){\circle*{3}}
\put(168,3){\circle*{3}} \put(100,15){\Large $Y$}

\end{picture}
\caption{Rc-graph $R$ and $Y$.} \label{ex2:R,T}
\end{figure}

Figures \ref{ex2:steps(3,0),(3,1),(2,0)}-\ref{ex2:final} show
rc-graphs $R(\ell,i)$ and tableaux of transposition $T(\ell,i)$.
Steps $(3,0)$, $(2,0)$ and $(1,0)$ are row-to-row steps. Steps
$(3,1)$, $(2,1)$, $(2,3)$ and $(1,2)$ are insertions. Steps
$(2,2)$, $(1,1)$, $(1,3)$ and $(1,4)$ are rectifications. We
circle all crossings of $R(\ell,i)$ with $i>0$, which are removed
or added by the current step. We also show by an arrow how
crossing move during rectifications.

Let us also recall that each row-to-row step $(\ell,0)$ constructs
the sequence of rc-graphs $S_j(\mathcal P(\ell+1))$ and then sets
$R(\ell,0)=S_0(\mathcal P(\ell+1))$. We omit the details of this
construction and refer to Section \ref{example1} for an example of
such construction. Also, after each row-to-row step
$w(R(\ell,0))=w(R_{\geq\ell})$, but rc-graphs $R(\ell,0)$ and
$R_{\geq\ell}$ could be different. For example, see step $(1,0)$
in Figure~\ref{ex2:steps(1,0),(1,1),(1,2)}.
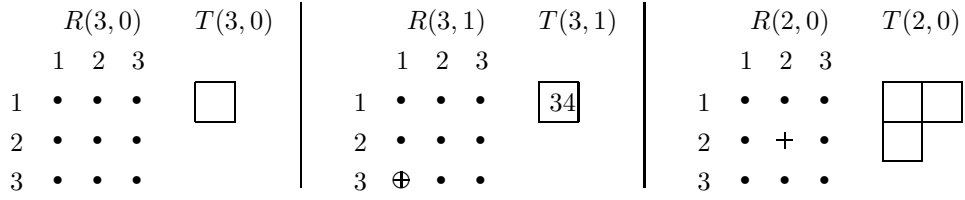
\begin{figure}[ht]
\begin{picture}(360,67)

\put (20,60){$R(3,0)$}


\put(0,29){1} \put(0,14){2} \put(0,-1){3}

\put(16,45){1} \put(31,45){2} \put(46,45){3}


\put(18,33){\circle*{3}} \put(33,33){\circle*{3}}
\put(33,18){\circle*{3}}
\put(33,3){\circle*{3}}\put(48,3){\circle*{3}}
\put(18,18){\circle*{3}} \put(18,3){\circle*{3}}
\put(48,18){\circle*{3}}\put(48,33){\circle*{3}}


\put (70,60){$T(3,0)$}

\put (70,40){\line(1,0){15}} \put (70,25){\line(1,0){15}}

\put (70,25){\line(0,1){15}} \put (85,25){\line(0,1){15}}

\put(110,0){\line(0,1){70}}


\put (150,60){$R(3,1)$}


\put(130,29){1} \put(130,14){2} \put(130,-1){3}

\put(146,45){1} \put(161,45){2} \put(176,45){3}


\put(148,33){\circle*{3}} \put(163,33){\circle*{3}}
\put(163,18){\circle*{3}}
\put(163,3){\circle*{3}}\put(178,3){\circle*{3}}
\put(148,18){\circle*{3}}
\put(178,18){\circle*{3}}\put(178,33){\circle*{3}}

\put (145,3){\line(1,0){6}} \put (148,0){\line(0,1){6}}

\put(148,3){\circle{6}}


\put (200,60){$T(3,1)$}

\put (200,40){\line(1,0){15}} \put (200,25){\line(1,0){15}}

\put (200,25){\line(0,1){15}} \put (215,25){\line(0,1){15}}


\put(204,29){$34$}


\put (280,60){$R(2,0)$}


\put(260,29){1} \put(260,14){2} \put(260,-1){3}

\put(276,45){1} \put(291,45){2} \put(306,45){3}


\put(278,33){\circle*{3}} \put(293,33){\circle*{3}}
\put(293,3){\circle*{3}}\put(308,3){\circle*{3}}
\put(278,18){\circle*{3}} \put(278,3){\circle*{3}}
\put(308,18){\circle*{3}}\put(308,33){\circle*{3}}


\put (290,18){\line(1,0){6}} \put (293,15){\line(0,1){6}}


\put (330,60){$T(2,0)$}

\put (330,40){\line(1,0){15}} \put (330,25){\line(1,0){15}} \put
(345,25){\line(1,0){15}} \put (345,40){\line(1,0){15}} \put
(330,10){\line(1,0){15}}

\put (330,25){\line(0,1){15}} \put (345,25){\line(0,1){15}} \put
(330,10){\line(0,1){15}} \put (345,10){\line(0,1){15}} \put
(360,25){\line(0,1){15}}
\put(240,0){\line(0,1){70}}

\end{picture}
\caption{Steps $(3,0)$, $(3,1)$ and $(2,0)$.}
\label{ex2:steps(3,0),(3,1),(2,0)}
\end{figure}

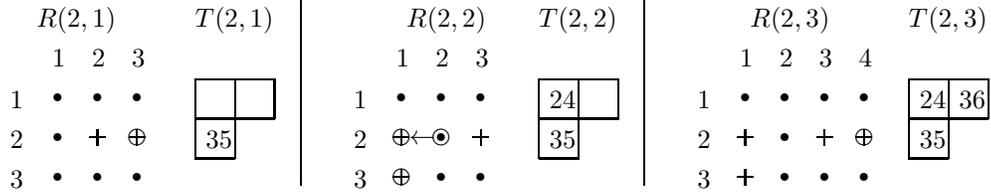
\begin{figure}[ht]
\begin{picture}(360,60)


\put (10,60){$R(2,1)$}


\put(0,29){1} \put(0,14){2} \put(0,-1){3}

\put(16,45){1} \put(31,45){2} \put(46,45){3}


\put(18,33){\circle*{3}} \put(33,33){\circle*{3}}
\put(18,3){\circle*{3}} \put(33,3){\circle*{3}}
\put(48,3){\circle*{3}} \put(18,18){\circle*{3}}
\put(48,33){\circle*{3}}

\put (30,18){\line(1,0){6}} \put (33,15){\line(0,1){6}}

\put (45,18){\line(1,0){6}} \put (48,15){\line(0,1){6}}

\put(48,18){\circle{6}}


\put (70,60){$T(2,1)$}

\put (70,40){\line(1,0){15}} \put (70,25){\line(1,0){15}} \put
(85,25){\line(1,0){15}} \put (85,40){\line(1,0){15}} \put
(70,10){\line(1,0){15}}

\put (70,25){\line(0,1){15}} \put (85,25){\line(0,1){15}} \put
(70,10){\line(0,1){15}} \put (85,10){\line(0,1){15}} \put
(100,25){\line(0,1){15}}


\put(74,14){$35$}

\put(110,0){\line(0,1){70}}


\put (150,60){$R(2,2)$}


\put(130,29){1} \put(130,14){2} \put(130,-1){3}

\put(146,45){1} \put(161,45){2} \put(176,45){3}


\put(148,33){\circle*{3}} \put(163,33){\circle*{3}}
\put(163,3){\circle*{3}}\put(178,3){\circle*{3}}
\put(163,18){\circle*{3}} \put(178,33){\circle*{3}}


\put (145,18){\line(1,0){6}} \put (148,15){\line(0,1){6}}

\put (175,18){\line(1,0){6}} \put (178,15){\line(0,1){6}}

\put (145,3){\line(1,0){6}} \put (148,0){\line(0,1){6}}

\put(148,18){\circle{6}} \put(163,18){\circle{6}}
\put(148,3){\circle{6}}

\put(151,15.5){$\leftarrow$}


\put (200,60){$T(2,2)$}

\put (200,40){\line(1,0){15}} \put (200,25){\line(1,0){15}} \put
(215,25){\line(1,0){15}} \put (215,40){\line(1,0){15}} \put
(200,10){\line(1,0){15}}

\put (200,25){\line(0,1){15}} \put (215,25){\line(0,1){15}} \put
(200,10){\line(0,1){15}} \put (215,10){\line(0,1){15}} \put
(230,25){\line(0,1){15}}

\put(204,14){$35$} \put(204,29){$24$}

\put(240,0){\line(0,1){70}}


\put (280,60){$R(2,3)$}


\put(260,29){1} \put(260,14){2} \put(260,-1){3}

\put(276,45){1} \put(291,45){2} \put(306,45){3} \put(321,45){4}


\put(278,33){\circle*{3}} \put(293,33){\circle*{3}}
\put(293,3){\circle*{3}} \put(308,3){\circle*{3}}
\put(293,18){\circle*{3}} \put(308,33){\circle*{3}}
\put(323,33){\circle*{3}} \put(323,3){\circle*{3}}

\put (275,18){\line(1,0){6}} \put (278,15){\line(0,1){6}}

\put (275,3){\line(1,0){6}} \put (278,0){\line(0,1){6}}

\put (305,18){\line(1,0){6}} \put (308,15){\line(0,1){6}}

\put (320,18){\line(1,0){6}} \put (323,15){\line(0,1){6}}

\put(323,18){\circle{6}}


\put (340,60){$T(2,3)$}

\put (340,40){\line(1,0){15}} \put (340,25){\line(1,0){15}} \put
(355,25){\line(1,0){15}} \put (355,40){\line(1,0){15}} \put
(340,10){\line(1,0){15}}

\put (340,25){\line(0,1){15}} \put (355,25){\line(0,1){15}} \put
(340,10){\line(0,1){15}} \put (355,10){\line(0,1){15}} \put
(370,25){\line(0,1){15}}


\put(344,14){$35$} \put(344,29){$24$} \put(359,29){$36$}

\end{picture}
\caption{Steps $(2,1)$, $(2,2)$ and $(2,3)$.}
\label{ex2:steps(2,1),(2,2),(2,3)}
\end{figure}

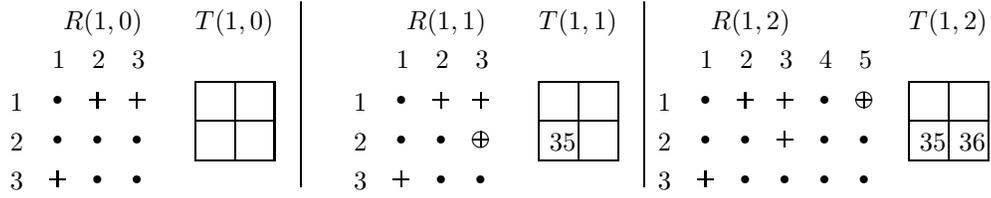
\begin{figure}[ht]
\begin{picture}(360,60)

\put (20,60){$R(1,0)$}


\put(0,29){1} \put(0,14){2} \put(0,-1){3}

\put(16,45){1} \put(31,45){2} \put(46,45){3}


\put(18,33){\circle*{3}} \put(33,3){\circle*{3}}
\put(48,3){\circle*{3}} \put(18,18){\circle*{3}}
\put(33,18){\circle*{3}} \put(48,18){\circle*{3}}


\put (15,3){\line(1,0){6}} \put (18,0){\line(0,1){6}}

\put (45,33){\line(1,0){6}} \put (48,30){\line(0,1){6}}

\put (30,33){\line(1,0){6}} \put (33,30){\line(0,1){6}}


\put (70,60){$T(1,0)$}

\put (70,40){\line(1,0){15}} \put (70,25){\line(1,0){15}} \put
(85,25){\line(1,0){15}} \put (85,40){\line(1,0){15}} \put
(70,10){\line(1,0){15}} \put (85,10){\line(1,0){15}}

\put (70,25){\line(0,1){15}} \put (85,25){\line(0,1){15}} \put
(70,10){\line(0,1){15}} \put (85,10){\line(0,1){15}} \put
(100,25){\line(0,1){15}} \put (100,10){\line(0,1){15}}
\put(110,0){\line(0,1){70}}


\put (150,60){$R(1,1)$}


\put(130,29){1} \put(130,14){2} \put(130,-1){3}

\put(146,45){1} \put(161,45){2} \put(176,45){3}


\put(148,33){\circle*{3}} \put(163,18){\circle*{3}}
\put(163,3){\circle*{3}} \put(178,3){\circle*{3}}
\put(148,18){\circle*{3}}

\put (145,3){\line(1,0){6}} \put (148,0){\line(0,1){6}}

\put (175,18){\line(1,0){6}} \put (178,15){\line(0,1){6}}

\put (160,33){\line(1,0){6}} \put (163,30){\line(0,1){6}}

\put (175,33){\line(1,0){6}} \put (178,30){\line(0,1){6}}

\put(178,18){\circle{6}}


\put (200,60){$T(1,1)$}

\put (200,40){\line(1,0){15}} \put (200,25){\line(1,0){15}} \put
(215,25){\line(1,0){15}} \put (215,40){\line(1,0){15}} \put
(200,10){\line(1,0){15}} \put (215,10){\line(1,0){15}}

\put (200,25){\line(0,1){15}} \put (215,25){\line(0,1){15}} \put
(200,10){\line(0,1){15}} \put (215,10){\line(0,1){15}} \put
(230,25){\line(0,1){15}} \put (230,10){\line(0,1){15}}


\put(204,14){$35$}


\put (265,60){$R(1,2)$}


\put(245,29){1} \put(245,14){2} \put(245,-1){3}

\put(261,45){1} \put(276,45){2} \put(291,45){3} \put(306,45){4}
\put(321,45){5}


\put(263,33){\circle*{3}} \put(293,3){\circle*{3}}
\put(263,18){\circle*{3}} \put(278,18){\circle*{3}}
\put(278,3){\circle*{3}} \put(308,3){\circle*{3}}
\put(308,18){\circle*{3}} \put(308,33){\circle*{3}}
\put(323,3){\circle*{3}} \put(323,18){\circle*{3}}


\put (260,3){\line(1,0){6}} \put (263,0){\line(0,1){6}}

\put (290,33){\line(1,0){6}} \put (293,30){\line(0,1){6}}

\put (275,33){\line(1,0){6}} \put (278,30){\line(0,1){6}}

\put (290,18){\line(1,0){6}} \put (293,15){\line(0,1){6}}

\put (290,33){\line(1,0){6}} \put (293,30){\line(0,1){6}}

\put (320,33){\line(1,0){6}} \put (323,30){\line(0,1){6}}

\put(323,33){\circle{6}}

\put(344,14){$35$} \put(359,14){$36$}

\put (340,60){$T(1,2)$}

\put (340,40){\line(1,0){15}} \put (340,25){\line(1,0){15}} \put
(355,25){\line(1,0){15}} \put (355,40){\line(1,0){15}} \put
(340,10){\line(1,0){15}} \put (355,10){\line(1,0){15}}

\put (340,25){\line(0,1){15}} \put (355,25){\line(0,1){15}} \put
(340,10){\line(0,1){15}} \put (355,10){\line(0,1){15}} \put
(370,25){\line(0,1){15}} \put (370,10){\line(0,1){15}}
\put(240,0){\line(0,1){70}}

\end{picture}
\caption{Steps $(1,0)$, $(1,1)$ and $(1,2)$.}
\label{ex2:steps(1,0),(1,1),(1,2)}
\end{figure}

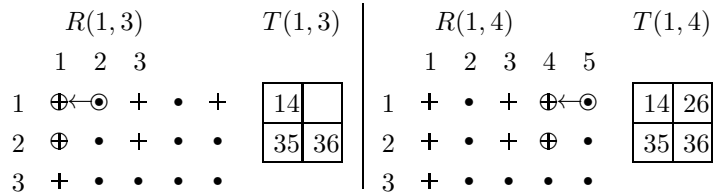
\begin{figure}[ht]
\begin{picture}(240,67)

\put(133,0){\line(0,1){70}}


\put (20,60){$R(1,3)$}


\put(0,29){1} \put(0,14){2} \put(0,-1){3}

\put(16,45){1} \put(31,45){2} \put(46,45){3}


\put(33,33){\circle*{3}} \put(33,18){\circle*{3}}
\put(33,3){\circle*{3}} \put(48,3){\circle*{3}}
 \put(63,3){\circle*{3}}
\put(63,18){\circle*{3}} \put(63,33){\circle*{3}}
\put(78,3){\circle*{3}} \put(78,18){\circle*{3}}

\put (15,3){\line(1,0){6}} \put (18,0){\line(0,1){6}}

\put (45,18){\line(1,0){6}} \put (48,15){\line(0,1){6}}

\put (15,33){\line(1,0){6}} \put (18,30){\line(0,1){6}}

\put (15,18){\line(1,0){6}} \put (18,15){\line(0,1){6}}

\put (45,33){\line(1,0){6}} \put (48,30){\line(0,1){6}}

\put (75,33){\line(1,0){6}} \put (78,30){\line(0,1){6}}

\put(18,18){\circle{6}}\put(18,33){\circle{6}}
\put(33,33){\circle{6}}

\put(21,30.5){$\leftarrow$}


\put (95,60){$T(1,3)$}

\put (95,40){\line(1,0){15}} \put (95,25){\line(1,0){15}} \put
(110,25){\line(1,0){15}} \put (110,40){\line(1,0){15}} \put
(95,10){\line(1,0){15}} \put (110,10){\line(1,0){15}}

\put (95,25){\line(0,1){15}} \put (110,25){\line(0,1){15}} \put
(95,10){\line(0,1){15}} \put (110,10){\line(0,1){15}} \put
(125,25){\line(0,1){15}} \put (125,10){\line(0,1){15}}


\put(99,14){$35$}\put(99,29){$14$}\put(114,14){$36$}

\put (160,60){$R(1,4)$}


\put(140,29){1} \put(140,14){2} \put(140,-1){3}

\put(156,45){1} \put(171,45){2} \put(186,45){3} \put(201,45){4}
\put(216,45){5}


\put(173,33){\circle*{3}} \put(173,18){\circle*{3}}
\put(173,3){\circle*{3}} \put(188,3){\circle*{3}}
 \put(203,3){\circle*{3}}
\put(218,33){\circle*{3}} \put(218,3){\circle*{3}}
\put(218,18){\circle*{3}}

\put (155,3){\line(1,0){6}} \put (158,0){\line(0,1){6}}

\put (185,18){\line(1,0){6}} \put (188,15){\line(0,1){6}}

\put (155,33){\line(1,0){6}} \put (158,30){\line(0,1){6}}

\put (155,18){\line(1,0){6}} \put (158,15){\line(0,1){6}}

\put (185,33){\line(1,0){6}} \put (188,30){\line(0,1){6}}

\put (200,33){\line(1,0){6}} \put (203,30){\line(0,1){6}}

\put (200,18){\line(1,0){6}} \put (203,15){\line(0,1){6}}

\put(203,18){\circle{6}}\put(203,33){\circle{6}}
\put(218,33){\circle{6}}

\put(206,30.5){$\leftarrow$}


\put (235,60){$T(1,4)$}

\put (235,40){\line(1,0){15}} \put (235,25){\line(1,0){15}} \put
(250,25){\line(1,0){15}} \put (250,40){\line(1,0){15}} \put
(235,10){\line(1,0){15}} \put (250,10){\line(1,0){15}}

\put (235,25){\line(0,1){15}} \put (250,25){\line(0,1){15}} \put
(235,10){\line(0,1){15}} \put (250,10){\line(0,1){15}} \put
(265,25){\line(0,1){15}} \put (265,10){\line(0,1){15}}


\put(239,14){$35$}\put(239,29){$14$}\put(254,14){$36$}
\put(254,29){$26$}

\end{picture}
\caption{Step $(1,3) $ and the final step $(1,4)$.}
\label{ex2:final}
\end{figure}

\newpage

\subsection{Example of insertion algorithm in case (\ref{case1})}
\label{example3} Let us now present an example in case shape of
$Y$ is a hook. Let $R$ and $Y$ be shown in Figure \ref{ex3:R,Y}.
In particular, $w(R)=(1,2,4,6,3,5,7,8,\dots)$ and
$w(Y)=(1,3,5,2,4,6,7,\dots)$, both are shuffles, but $w(R)$ has
descent at $4$, while $w(Y)$ has descent at $3$, so case
(\ref{case2}) does not apply.

\begin{figure}[ht]
\begin{picture}(200,56)



\put(20,29){1} \put(20,14){2} \put(20,-1){3}

\put(36,45){1} \put(51,45){2} \put(66,45){3}


\put (50,3){\line(1,0){6}} \put (53,0){\line(0,1){6}}

\put (65,33){\line(1,0){6}} \put (68,30){\line(0,1){6}}

\put (65,3){\line(1,0){6}} \put (68,0){\line(0,1){6}}


\put(38,33){\circle*{3}} \put(53,18){\circle*{3}}
\put(53,33){\circle*{3}} \put(38,3){\circle*{3}}
\put(68,18){\circle*{3}} \put(38,18){\circle*{3}}

\put(00,15){\Large $R$}



\put(120,29){1} \put(120,14){2} \put(120,-1){3}

\put(136,45){1} \put(151,45){2} \put(166,45){3} \put(181,45){4}


\put (135,3){\line(1,0){6}} \put (138,0){\line(0,1){6}}

\put (150,33){\line(1,0){6}} \put (153,30){\line(0,1){6}}

\put (180,33){\line(1,0){6}} \put (183,30){\line(0,1){6}}


\put(138,33){\circle*{3}} \put(153,18){\circle*{3}}
\put(168,33){\circle*{3}} \put(153,3){\circle*{3}}
\put(168,3){\circle*{3}} \put(138,18){\circle*{3}}
\put(183,3){\circle*{3}} \put(168,18){\circle*{3}}
\put(183,18){\circle*{3}}

\put(100,15){\Large $Y$}

\end{picture}
\caption{Rc-graph $R$ and $Y$.} \label{ex3:R,Y}
\end{figure}
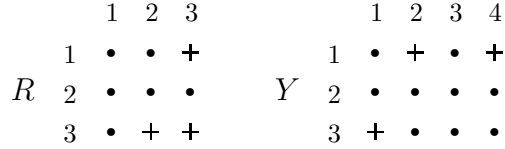

Figures \ref{ex3:steps(3,0),(3,1),(2,0)}-\ref{ex3:final} contain
the results of all the steps of the algorithm.  Steps $(3,0)$,
$(2,0)$ and $(1,0)$ are row-to-row steps, steps $(3,1)$, $(1,1)$
and $(1,3)$ are insertions, while $(2,1)$ and $(1,2)$ are
rectifications. Notice that step $(1,2)$ is the only step, where
the second situation of Figure \ref{remove-crossing} occurs.

\begin{figure}[ht]
\begin{picture}(360,70)

\put (20,60){$R(3,0)$}


\put(0,29){1} \put(0,14){2} \put(0,-1){3}

\put(16,45){1} \put(31,45){2} \put(46,45){3}


\put(18,33){\circle*{3}} \put(33,33){\circle*{3}}
\put(33,18){\circle*{3}} \put(18,3){\circle*{3}}
\put(18,18){\circle*{3}}
\put(48,18){\circle*{3}}\put(48,33){\circle*{3}}


\put (30,3){\line(1,0){6}} \put (33,0){\line(0,1){6}}

\put (45,3){\line(1,0){6}} \put (48,0){\line(0,1){6}}


\put (70,60){$T(3,0)$}

\put (70,40){\line(1,0){15}} \put (70,25){\line(1,0){15}}

\put (70,25){\line(0,1){15}} \put (85,25){\line(0,1){15}}

\put(110,0){\line(0,1){70}}


\put (150,60){$R(3,1)$}


\put(130,29){1} \put(130,14){2} \put(130,-1){3}

\put(146,45){1} \put(161,45){2} \put(176,45){3}


\put(148,33){\circle*{3}} \put(163,33){\circle*{3}}
\put(163,18){\circle*{3}} \put(148,18){\circle*{3}}
\put(178,18){\circle*{3}}\put(178,33){\circle*{3}}

\put (145,3){\line(1,0){6}} \put (148,0){\line(0,1){6}}

\put (160,3){\line(1,0){6}} \put (163,0){\line(0,1){6}}

\put (175,3){\line(1,0){6}} \put (178,0){\line(0,1){6}}

\put(148,3){\circle{6}}


\put (200,60){$T(3,1)$}

\put (200,40){\line(1,0){15}} \put (200,25){\line(1,0){15}}

\put (200,25){\line(0,1){15}} \put (215,25){\line(0,1){15}}


\put(204,29){$34$}


\put (280,60){$R(2,0)$}


\put(260,29){1} \put(260,14){2} \put(260,-1){3}

\put(276,45){1} \put(291,45){2} \put(306,45){3}


\put(278,33){\circle*{3}} \put(293,33){\circle*{3}}
\put(293,18){\circle*{3}} \put(278,3){\circle*{3}}
\put(278,18){\circle*{3}}
\put(308,18){\circle*{3}}\put(308,33){\circle*{3}}


\put (290,3){\line(1,0){6}} \put (293,0){\line(0,1){6}}

\put (305,3){\line(1,0){6}} \put (308,0){\line(0,1){6}}


\put (330,60){$T(2,0)$}

\put (330,40){\line(1,0){15}} \put (330,25){\line(1,0){15}}

\put (330,25){\line(0,1){15}} \put (345,25){\line(0,1){15}}

\put(240,0){\line(0,1){70}}

\end{picture}
\caption{Steps $(3,0)$, $(3,1)$ and $(2,0)$.}
\label{ex3:steps(3,0),(3,1),(2,0)}
\end{figure}
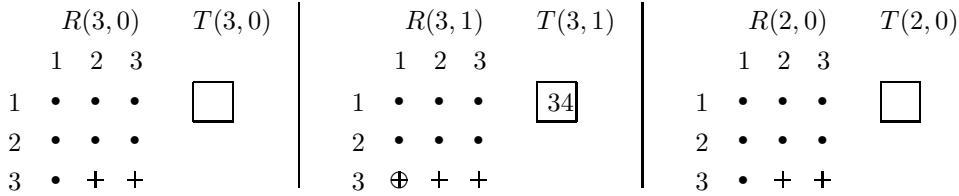

\begin{figure}[ht]
\begin{picture}(360,70)

\put(110,0){\line(0,1){70}}


\put (20,60){$R(2,1)$}


\put(0,29){1} \put(0,14){2} \put(0,-1){3}

\put(16,45){1} \put(31,45){2} \put(46,45){3}


\put(18,33){\circle*{3}} \put(33,33){\circle*{3}}
\put(33,18){\circle*{3}} \put(18,18){\circle*{3}}
\put(48,18){\circle*{3}}\put(48,33){\circle*{3}}

\put (15,3){\line(1,0){6}} \put (18,0){\line(0,1){6}}

\put (30,3){\line(1,0){6}} \put (33,0){\line(0,1){6}}

\put (45,3){\line(1,0){6}} \put (48,0){\line(0,1){6}}

\put(18,3){\circle{6}}


\put (70,60){$T(2,1)$}

\put (70,40){\line(1,0){15}} \put (70,25){\line(1,0){15}}

\put (70,25){\line(0,1){15}} \put (85,25){\line(0,1){15}}


\put(74,29){$34$}

\put (150,60){$R(1,0)$}


\put(130,29){1} \put(130,14){2} \put(130,-1){3}

\put(16,45){1} \put(31,45){2} \put(46,45){3}


\put(148,33){\circle*{3}} \put(163,33){\circle*{3}}
\put(163,18){\circle*{3}} \put(148,3){\circle*{3}}
\put(148,18){\circle*{3}} \put(178,18){\circle*{3}}


\put (160,3){\line(1,0){6}} \put (163,0){\line(0,1){6}}

\put (175,3){\line(1,0){6}} \put (178,0){\line(0,1){6}}

\put (175,33){\line(1,0){6}} \put (178,30){\line(0,1){6}}


\put (200,60){$T(1,0)$}

\put (200,40){\line(1,0){15}} \put (200,25){\line(1,0){15}} \put
(215,25){\line(1,0){15}} \put (215,40){\line(1,0){15}} \put
(200,10){\line(1,0){15}}

\put (200,25){\line(0,1){15}} \put (215,25){\line(0,1){15}} \put
(200,10){\line(0,1){15}} \put (215,10){\line(0,1){15}} \put
(230,25){\line(0,1){15}}
\put(240,0){\line(0,1){70}}


\put (280,60){$R(1,1)$}


\put(260,29){1} \put(260,14){2} \put(260,-1){3}

\put(276,45){1} \put(291,45){2} \put(306,45){3} \put(321,45){4}


\put(278,33){\circle*{3}} \put(293,33){\circle*{3}}
\put(293,18){\circle*{3}} \put(278,18){\circle*{3}}
\put(308,18){\circle*{3}}\put(278,3){\circle*{3}}
\put(323,18){\circle*{3}}\put(323,3){\circle*{3}}

\put (305,33){\line(1,0){6}} \put (308,30){\line(0,1){6}}

\put (320,33){\line(1,0){6}} \put (323,30){\line(0,1){6}}

\put (290,3){\line(1,0){6}} \put (293,0){\line(0,1){6}}

\put (305,3){\line(1,0){6}} \put (308,0){\line(0,1){6}}

\put(323,33){\circle{6}}


\put (340,60){$T(1,1)$}

\put (340,40){\line(1,0){15}} \put (340,25){\line(1,0){15}} \put
(355,25){\line(1,0){15}} \put (355,40){\line(1,0){15}} \put
(340,10){\line(1,0){15}}

\put (340,25){\line(0,1){15}} \put (355,25){\line(0,1){15}} \put
(340,10){\line(0,1){15}} \put (355,10){\line(0,1){15}} \put
(370,25){\line(0,1){15}}


\put(344,14){$36$}

\end{picture}
\caption{Steps $(2,1)$, $(1,0)$ and $(1,1)$.}
\label{ex3:steps(2,1),(1,0),(1,1)}
\end{figure}
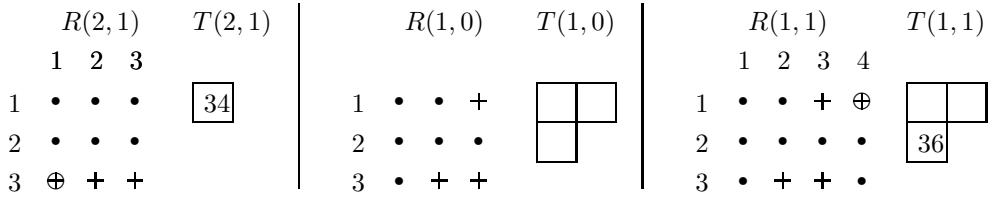

\begin{figure}[ht]
\begin{picture}(240,70)

\put (20,60){$R(1,2)$}


\put(0,29){1} \put(0,14){2} \put(0,-1){3}

\put(16,45){1} \put(31,45){2} \put(46,45){3} \put(61,45){4}


\put(18,33){\circle*{3}} \put(63,33){\circle*{3}}
\put(33,18){\circle*{3}} \put(63,3){\circle*{3}}
\put(18,18){\circle*{3}} \put(48,18){\circle*{3}}
\put(63,18){\circle*{3}}


\put (30,3){\line(1,0){6}} \put (33,0){\line(0,1){6}}

\put (15,3){\line(1,0){6}} \put (18,0){\line(0,1){6}}

\put (30,33){\line(1,0){6}} \put (33,30){\line(0,1){6}}

\put (45,3){\line(1,0){6}} \put (48,0){\line(0,1){6}}

\put (45,33){\line(1,0){6}} \put (48,30){\line(0,1){6}}

\put(63,33){\circle{6}} \put(33,33){\circle{6}}
\put(18,3){\circle{6}}

\put (80,60){$T(1,2)$}

\put (80,40){\line(1,0){15}} \put (80,25){\line(1,0){15}} \put
(95,25){\line(1,0){15}} \put (95,40){\line(1,0){15}} \put
(80,10){\line(1,0){15}}

\put (80,25){\line(0,1){15}} \put (95,25){\line(0,1){15}} \put
(80,10){\line(0,1){15}} \put (95,10){\line(0,1){15}} \put
(110,25){\line(0,1){15}}

\put(84,14){$34$} \put(84,29){$24$}

\put(35,30.5){$\longleftarrow$}

\put(45,33){\line(1,0){15}}

\put(120,0){\line(0,1){70}}


\put (160,60){$R(1,3)$}


\put(140,29){1} \put(140,14){2} \put(140,-1){3}

\put(156,45){1} \put(171,45){2} \put(186,45){3} \put(201,45){4}
\put(216,45){5} \put(231,45){6}

\put(173,18){\circle*{3}} \put(233,3){\circle*{3}}
\put(203,3){\circle*{3}} \put(203,18){\circle*{3}}
\put(203,33){\circle*{3}} \put(218,3){\circle*{3}}
\put(218,18){\circle*{3}} \put(158,18){\circle*{3}}
\put(158,33){\circle*{3}} \put(188,18){\circle*{3}}
\put(233,18){\circle*{3}} \put(218,33){\circle*{3}}
\put (155,3){\line(1,0){6}} \put (158,0){\line(0,1){6}}

\put (170,3){\line(1,0){6}} \put (173,0){\line(0,1){6}}

\put (185,3){\line(1,0){6}} \put (188,0){\line(0,1){6}}

\put (230,33){\line(1,0){6}} \put (233,30){\line(0,1){6}}

\put (170,33){\line(1,0){6}} \put (173,30){\line(0,1){6}}

\put (185,33){\line(1,0){6}} \put (188,30){\line(0,1){6}}

\put(233,33){\circle{6}}

\put (245,60){$T(1,3)$}

\put (245,40){\line(1,0){15}} \put (245,25){\line(1,0){15}} \put
(260,25){\line(1,0){15}} \put (260,40){\line(1,0){15}} \put
(245,10){\line(1,0){15}}

\put (245,25){\line(0,1){15}} \put (260,25){\line(0,1){15}} \put
(245,10){\line(0,1){15}} \put (260,10){\line(0,1){15}} \put
(275,25){\line(0,1){15}}

\put(249,14){$34$} \put(249,29){$24$} \put(264,29){$37$}

\end{picture}
\caption{Step $(1,2)$ and the final step $(1,3)$.}
\label{ex3:final}
\end{figure}
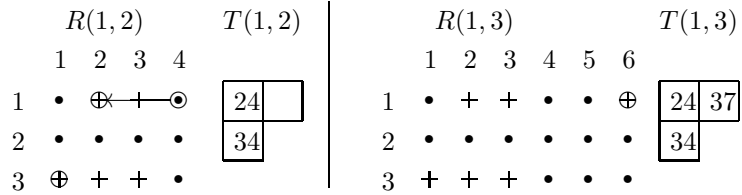

\newpage

\subsection{Another example in case (\ref{case1})}
The last example is for rc-graphs $R$ and $Y$ defined in Figure
\ref{ex4:R,Y}. In this case $w(R)=(1,2,5,4,6,3,7,8,\dots)$ and
$w(Y)=(1,4,2,3,5,6,\dots)=v((2,0),2)$, the shape of $Y$ is a row,
while $w(R)$ is a permutation with two descents.

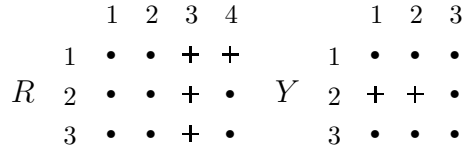
\begin{figure}[ht]
\begin{picture}(200,56)



\put(20,29){1} \put(20,14){2} \put(20,-1){3}

\put(36,45){1} \put(51,45){2} \put(66,45){3} \put(81,45){4}


\put (65,18){\line(1,0){6}} \put (68,15){\line(0,1){6}}

\put (65,33){\line(1,0){6}} \put (68,30){\line(0,1){6}}

\put (65,3){\line(1,0){6}} \put (68,0){\line(0,1){6}}

\put (80,33){\line(1,0){6}} \put (83,30){\line(0,1){6}}


\put(38,33){\circle*{3}} \put(53,18){\circle*{3}}
\put(53,33){\circle*{3}} \put(38,3){\circle*{3}}
\put(53,3){\circle*{3}} \put(38,18){\circle*{3}}
\put(83,3){\circle*{3}}\put(83,18){\circle*{3}}

\put(00,15){\Large $R$}



\put(120,29){1} \put(120,14){2} \put(120,-1){3}

\put(136,45){1} \put(151,45){2} \put(166,45){3}


\put (135,18){\line(1,0){6}} \put (138,15){\line(0,1){6}}

\put (150,18){\line(1,0){6}} \put (153,15){\line(0,1){6}}


\put(138,33){\circle*{3}} \put(153,33){\circle*{3}}
\put(168,33){\circle*{3}} \put(153,3){\circle*{3}}
\put(168,3){\circle*{3}} \put(138,3){\circle*{3}}
\put(168,18){\circle*{3}}

\put(100,15){\Large $Y$}

\end{picture}
\caption{Rc-graph $R$ and $Y$.} \label{ex4:R,Y}
\end{figure}
The steps of the algorithm are shown in Figures
\ref{ex4:steps(2,0),(2,1),(2,2)}-\ref{ex4:final}. Steps $(2,0)$
and $(1,0)$ are row-to-row steps, steps $(2,1)$ and $(2,1)$ are
insertions, while steps $(1,1)$ and $(1,2)$ are rectifications.
Notice that step $(1,2)$ is the only step where the second case of
Figure \ref{allowed-insertions-rect} occurs.

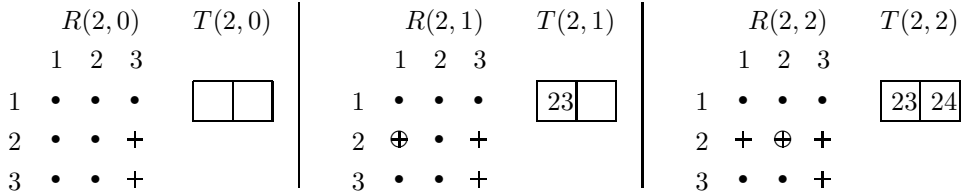
\begin{figure}[ht]
\begin{picture}(360,70)

\put (20,60){$R(2,0)$}


\put(0,29){1} \put(0,14){2} \put(0,-1){3}

\put(16,45){1} \put(31,45){2} \put(46,45){3}


\put(18,33){\circle*{3}} \put(33,33){\circle*{3}}
\put(33,18){\circle*{3}} \put(18,3){\circle*{3}}
\put(18,18){\circle*{3}}
\put(33,3){\circle*{3}}\put(48,33){\circle*{3}}


\put (45,18){\line(1,0){6}} \put (48,15){\line(0,1){6}}

\put (45,3){\line(1,0){6}} \put (48,0){\line(0,1){6}}


\put (70,60){$T(2,0)$}

\put (70,40){\line(1,0){15}} \put (70,25){\line(1,0){15}} \put
(85,40){\line(1,0){15}} \put (85,25){\line(1,0){15}}

\put (70,25){\line(0,1){15}} \put (85,25){\line(0,1){15}} \put
(100,25){\line(0,1){15}}

\put(110,0){\line(0,1){70}}


\put (150,60){$R(2,1)$}


\put(130,29){1} \put(130,14){2} \put(130,-1){3}

\put(146,45){1} \put(161,45){2} \put(176,45){3}


\put(148,33){\circle*{3}} \put(163,33){\circle*{3}}
\put(163,18){\circle*{3}} \put(148,3){\circle*{3}}
\put(163,3){\circle*{3}}\put(178,33){\circle*{3}}

\put (145,18){\line(1,0){6}} \put (148,15){\line(0,1){6}}

\put (175,18){\line(1,0){6}} \put (178,15){\line(0,1){6}}

\put (175,3){\line(1,0){6}} \put (178,0){\line(0,1){6}}

\put(148,18){\circle{6}}


\put (200,60){$T(2,1)$}

\put (200,40){\line(1,0){15}} \put (200,25){\line(1,0){15}} \put
(215,40){\line(1,0){15}} \put (215,25){\line(1,0){15}}

\put (200,25){\line(0,1){15}} \put (215,25){\line(0,1){15}} \put
(230,25){\line(0,1){15}}


\put(204,29){$23$}


\put (280,60){$R(2,2)$}


\put(260,29){1} \put(260,14){2} \put(260,-1){3}

\put(276,45){1} \put(291,45){2} \put(306,45){3}


\put(278,33){\circle*{3}} \put(293,33){\circle*{3}}
\put(293,3){\circle*{3}} \put(278,3){\circle*{3}}
\put(308,33){\circle*{3}}


\put (275,18){\line(1,0){6}} \put (278,15){\line(0,1){6}}

\put (290,18){\line(1,0){6}} \put (293,15){\line(0,1){6}}

\put(305,18){\line(1,0){6}} \put (308,15){\line(0,1){6}}

\put (305,3){\line(1,0){6}} \put (308,0){\line(0,1){6}}

\put(293,18){\circle{6}}


\put (330,60){$T(2,2)$}

\put (330,40){\line(1,0){15}} \put (330,25){\line(1,0){15}} \put
(345,40){\line(1,0){15}} \put (345,25){\line(1,0){15}}

\put (330,25){\line(0,1){15}} \put (345,25){\line(0,1){15}} \put
(360,25){\line(0,1){15}}

\put(334,29){$23$} \put(349,29){$24$}

\put(240,0){\line(0,1){70}}

\end{picture}
\caption{Steps $(2,0)$, $(2,1)$ and $(2,2)$.}
\label{ex4:steps(2,0),(2,1),(2,2)}
\end{figure}

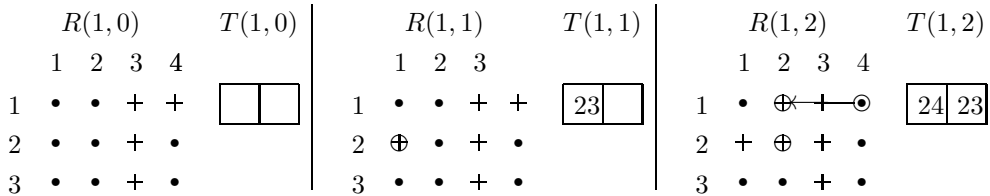
\begin{figure}[ht]
\begin{picture}(360,70)

\put(115,0){\line(0,1){70}}


\put (20,60){$R(1,0)$}


\put(0,29){1} \put(0,14){2} \put(0,-1){3}

\put(16,45){1} \put(31,45){2} \put(46,45){3} \put(61,45){4}

\put(18,33){\circle*{3}} \put(33,33){\circle*{3}}
\put(33,18){\circle*{3}} \put(18,18){\circle*{3}}
\put(33,3){\circle*{3}} \put(18,3){\circle*{3}}
\put(63,3){\circle*{3}} \put(63,18){\circle*{3}}

\put (45,18){\line(1,0){6}} \put (48,15){\line(0,1){6}}

\put (45,3){\line(1,0){6}} \put (48,0){\line(0,1){6}}

\put (45,33){\line(1,0){6}} \put (48,30){\line(0,1){6}}

\put (60,33){\line(1,0){6}} \put (63,30){\line(0,1){6}}


\put (80,60){$T(1,0)$}

\put (80,40){\line(1,0){15}} \put (80,25){\line(1,0){15}} \put
(95,40){\line(1,0){15}} \put (95,25){\line(1,0){15}}

\put (80,25){\line(0,1){15}} \put (95,25){\line(0,1){15}} \put
(110,25){\line(0,1){15}}

\put (150,60){$R(1,1)$}


\put(130,29){1} \put(130,14){2} \put(130,-1){3}

\put(146,45){1} \put(161,45){2} \put(176,45){3} \put(61,45){4}


\put(148,33){\circle*{3}} \put(163,33){\circle*{3}}
\put(163,18){\circle*{3}} \put(148,3){\circle*{3}}
\put(163,3){\circle*{3}} \put(193,18){\circle*{3}}
\put(193,3){\circle*{3}}

\put (145,18){\line(1,0){6}} \put (148,15){\line(0,1){6}}

\put (175,18){\line(1,0){6}} \put (178,15){\line(0,1){6}}

\put (175,3){\line(1,0){6}} \put (178,0){\line(0,1){6}}

\put (175,33){\line(1,0){6}} \put (178,30){\line(0,1){6}}

\put (190,33){\line(1,0){6}} \put (193,30){\line(0,1){6}}

\put(148,18){\circle{6}}


\put (210,60){$T(1,1)$}

\put (210,40){\line(1,0){15}} \put (210,25){\line(1,0){15}} \put
(225,25){\line(1,0){15}} \put (225,40){\line(1,0){15}}

\put (210,25){\line(0,1){15}} \put (225,25){\line(0,1){15}}  \put
(240,25){\line(0,1){15}}

\put(214,29){$23$}
\put(245,0){\line(0,1){70}}


\put (280,60){$R(1,2)$}


\put(260,29){1} \put(260,14){2} \put(260,-1){3}

\put(276,45){1} \put(291,45){2} \put(306,45){3} \put(321,45){4}


\put(278,33){\circle*{3}} \put(323,33){\circle*{3}}
\put(293,3){\circle*{3}} \put(278,3){\circle*{3}}
\put(323,18){\circle*{3}}\put(323,3){\circle*{3}}

\put (305,33){\line(1,0){6}} \put (308,30){\line(0,1){6}}

\put (305,18){\line(1,0){6}} \put (308,15){\line(0,1){6}}

\put (275,18){\line(1,0){6}} \put (278,15){\line(0,1){6}}

\put (290,18){\line(1,0){6}} \put (293,15){\line(0,1){6}}

\put (290,33){\line(1,0){6}} \put (293,30){\line(0,1){6}}

\put (305,3){\line(1,0){6}} \put (308,0){\line(0,1){6}}

\put(323,33){\circle{6}} \put(293,33){\circle{6}}
\put(293,18){\circle{6}}


\put (340,60){$T(1,2)$}

\put (340,40){\line(1,0){15}} \put (340,25){\line(1,0){15}} \put
(355,25){\line(1,0){15}} \put (355,40){\line(1,0){15}}

\put (340,25){\line(0,1){15}} \put (355,25){\line(0,1){15}}  \put
(370,25){\line(0,1){15}}


\put(359,29){$23$}\put(344,29){$24$}

\put(295,30.5){$\longleftarrow$}

\put(305,33){\line(1,0){15}}

\end{picture}
\caption{Steps $(1,0)$, $(1,1)$ and final step $(1,2)$.}
\label{ex4:final}
\end{figure}

\end{document}